\documentclass[a4paper,12pt]{article}
\usepackage[a4paper, total={15cm, 23.7cm}]{geometry}
\usepackage{amssymb,amsthm,amsmath,amsfonts,color,tikz}
\usepackage{comment}
\usepackage[normalem]{ulem}

\usepackage{graphicx}

\usepackage{hyperref}

\usepackage{lineno}

\usepackage[T1]{fontenc}
\usepackage{tikz}
\usetikzlibrary{calc,matrix}

\newtheorem{theorem}{Theorem}[section]
\newtheorem{lemma}[theorem]{Lemma}
\newtheorem{claim}[theorem]{Claim}
\newtheorem{proposition}[theorem]{Proposition}
\newtheorem{corollary}[theorem]{Corollary}

\newtheorem{definition}{Definition}[section]
\newtheorem{remark} {Remark}[section]

\newcommand{\bp}{\begin{proof}}
\newcommand{\ep}{\end{proof}}

\def\vs{\vskip4pt}

\begin{document}

\begin{sloppypar}

\allowdisplaybreaks

\title{Multiple existence and qualitative property of nodal solutions for coupled elliptic equations}

\author{{Haoyu Li} \\
Departamento de Matem\'atica\\
Universidade Federal de S\~ao Carlos\\
S\~ao Carlos-SP, 13565-905, Brazil\\
\texttt{hyli1994@hotmail.com}\\
\\
Zhi-Qiang Wang\\
Department of Mathematics and Statistics\\
Utah State University\\
Logan, Utah 84322, USA\\
\texttt{zhi-qiang.wang@usu.edu}}
\maketitle
\vskip1truecm\indent

\begin{abstract}\noindent
The paper studies nodal solutions having prescribed componentwise nodal data for the following coupled nonlinear elliptic equations
\begin{equation}
    \left\{
   \begin{array}{lr}
   -{\Delta}u_{j}+ u_{j}= u^{3}_{j}+\beta\sum_{i=1, i\neq j}^N u_{j}u_{i}^{2} \,\,\,\,\,\,\,  \mbox{in}\ \Omega,\nonumber\\
     u_{j}\in H_{0,r}^{1}(\Omega), \,\,\,\,\,\,\,\,j=1,\dots,N.\nonumber
   \end{array}
   \right.
\end{equation}
Here, $\Omega\subset\mathbb{R}^n$ is a bounded and radial domain with $n=2,3$. The coupling constant $\beta\leq-1$ is in the repulsive regime.
We investigate the solution structure for both positive and nodal solutions, proving multiple existence of solutions with prescribed nodal data and providing qualitative estimates for the nodal numbers of the inter-componentwise differences of solutions with both upper and lower bounds. Our general framework is for nodal solutions though our results are new also for positive solutions.
\end{abstract}


\tableofcontents

\section{Introduction}
\subsection{An overview}
In this paper, we consider the following nonlinear elliptic problem
\begin{equation}\label{e:AAAA}
    \left\{
   \begin{array}{lr}
   -{\Delta}u_{j}+ u_{j}= u^{3}_{j}+\beta\sum_{i=1, i\neq j}^N u_{j}u_{i}^{2} \,\,\,\,\,\,\,  \mbox{in}\ \Omega,\\
     u_{j}\in H_{0}^{1}(\Omega), \,\,\,\,\,\,\,\,j=1,\dots,N.
   \end{array}
   \right.
\end{equation}
Here, $\Omega\subset\mathbb{R}^n$ is a bounded and radial domain (e.g., a ball or annulus) with $n=2,3$, and the coupling constant $\beta\leq-1$. It is also referred as the coupled Schr\"odinger equations since it is the system for standing waves of the time dependent nonlinear Schr\"odinger equations:
\begin{equation}
    \left\{
   \begin{array}{lr}
     i\partial_{t}\Phi_{j}+{\Delta}\Phi_{j}+|\Phi_{j}|^2 \Phi_{j}+\beta\sum_{i\neq j}\Phi_{j}|\Phi_{i}|^{2}=0 \,\,\,\,\,\,\,  \mbox{in}\ \Omega ,\\
     \Phi_{j}(t,x)\in\mathbb{C},\Phi_{j}(t,x)|_{x\in\partial\Omega}=0 \,\,\,\,\,\,\,\,j=1,\dots,N.
   \end{array}
   \right.
\end{equation}
Such kind of systems arise from the Bose-Einstein condensate, cf. \cite{AkhmedievAnkiewicz1999,MitchellSegev1997}. Problem (\ref{e:AAAA}) is regarded as attractive if $\beta>0$, repulsive if $\beta<0$. In this paper, we focus on the repulsive regime.

In recent decades Problem (\ref{e:AAAA}) has been studied extensively, following with the seminar work of Lin-Wei \cite{LinWei2005}. It seems impossible to cover all the references and we discuss here some relevant ones as motivations to our work. The solvability of Problem (\ref{e:AAAA}) is aided with the help of the mountain-pass theorem, cf. \cite{RabinBook}. Nevertheless, it is know that (e.g., \cite{BartschWang2006}) the mountain pass solution could be semi-trivial, i.e., some of the components in the solution are zero while others are not. With the help of the invariant set of descending flow and the Nehari method, Liu and Wang \cite{LiuWang2008,LiuWang2010} overcame this difficulty and obtained an unbounded sequence of non-trivial solutions, i.e., solutions without any zero component. A significant difference between Problem (\ref{e:AAAA}) and its scalar field counterpart $-\Delta u+u=|u|^{p-2}u$ arises in \cite{WeiWeth2008,DancerWeiWeth2010,BartschDancerWang2010,TianWang2011}, where the authors found an infinite sequence of positive solutions under suitable assumptions.
The a priori estimate \cite{GidasSpruck1981} and the uniqueness result \cite{GidasNiNirenberg1979,Kwong1989} rule out the possibility of multiplicity of positive solutions for the scalar field equation.
This indicates that the full structure of solutions to Problem (\ref{e:AAAA}) when we also consider sign-changing solutions is far more complicated than that of their scalar counterpart.

For nodal solutions, several results about Problem (\ref{e:AAAA}) have appeared. The works in \cite{SatoWang2013, LiuLiuWang2015, ChenLinZou1016} found infinitely many nodal solutions or solutions with some of the components positive and the rest nodal. In the radial setting, solutions with a prescribed number of nodes for each component can be found in \cite{LiuWang2019,LiWang2021}. Furthermore, the authors of the current paper proved in \cite{LiWang2021} that with each set of component-wisely prescribed nodal numbers there exist an infinite sequence of solutions carrying the same set of nodal data. The approach in \cite{LiWang2021} is based on a special symmetric mountain-pass procedure by using a parabolic flow serving as a descending flow of the variational formulation.

In this paper we continue our investigation on the solution structure of both positive and nodal solutions in the radial setting,
We aim to achieve the following several goals. We will give a more detailed classification and properties of both positive and nodal solutions, in terms of nodal numbers for the components of solutions and for the differences between components of solutions.
We provide a new approach which enable us to strengthen and extend the existence results of nodal and positive solutions under weaker conditions.
Furthermore, our results are new even for positive solutions.

To demonstrate the spirit of our results we consider a special case here first, the full results will be presented in Section 1.2.
Let us consider the following coupled system of three equations, we write $(u_1,u_2,u_3)$ as $(u,v,w)$,
\begin{equation}\label{e:3Coupled}
    \left\{
   \begin{array}{lr}
     -{\Delta}u+u=u^{3}+\beta uv^2+\beta uw^2 \,\,\,  \mbox{in}\ B_1,\\
     -\Delta v+v =v^3 +\beta u^2 v+\beta vw^2 \,\,\,  \mbox{in}\ B_1,\\
     -\Delta w+w =w^3 +\beta u^2 w+\beta v^2w \,\,\,  \mbox{in}\ B_1,\\
     u,v,w\in H_{0,r}^{1}(B_1).
   \end{array}
   \right.
\end{equation}
Here, $\beta\leq-1$ and $B_1$ is the unit ball in $\mathbb R^n$ with $n=2,3$.
For a continuous radial function $u$, $n(u)$ denotes its number of zeros whenever it makes sense.
Our general work yields the following result for Problem (\ref{e:3Coupled}).
\begin{theorem}\label{t:000}
There exist four unbounded sequences of positive radial solutions to Problem (\ref{e:3Coupled}): $\{(u^l_s, v^l_s, w^l_s)\}_{s=1}^\infty$ for $l = 1, 2, 3, 4$, such that the following distinctly different inter-componentwise nodal properties hold:
\begin{itemize}
  \item [$(1).$] $\lim_{s\to\infty} n(u^1_s - v^1_s)=\infty$,
  and
  $ n(u^1_s - w^1_s)=n(v^1_s - w^1_s)=1$ for any $s\geq 1$;
  \item [$(2).$] $\lim_{s\to\infty}n(u^2_s - w^2_s)=\infty$,
  and $n(u^2_s - v^2_s)=n(v^2_s - w^2_s)=1$ for any $s\geq 1$;
  \item [$(3).$] $\lim_{s\to\infty} n(v^3_s - w^3_s)=\infty$,
  and $n(u^3_s - w^3_s)=n(u^3_s - v^3_s)=1$ any $s\geq 1$;
  \item [$(4).$] $n(u^4_s - v^4_s)$, $n(u^4_s-w^4_s)$ and $n(v^4_s - w^4_s)$ all tend to infinity as $s\to\infty$.
\end{itemize}
\end{theorem}
We remark that the importance of comparison of nodal numbers between components is observed in \cite{QuittnerNew}
for studying qualitative property of solutions, where Quittner obtained a priori estimate with prescribed nodal numbers for both components and comparisons.
The above result gives distinctly different infinite sequences of positive solutions, therefore further classifying the structure of positive solutions. The classification is based on the comparison of nodal numbers for the differences between the components of solutions.

\subsection{Main results and main tools}
We now return to Problem (\ref{e:AAAA}) and give the full results of the paper.
Throughout this article, we will denote the whole vector-valued functions by the capital letters such as $U$, $V$ or $W$. The energy functional of the problem (\ref{e:AAAA}) is defined as
\begin{align}\label{energy}
I(U)=\frac{1}{2}\sum_{j=1}^N\int|\nabla u_j|^2+ u^2_j-\frac{1}{4}\sum_{j=1}^N\int u_j^4+\beta\sum_{i\neq j}u_i^2 u_j^2
\end{align}
with $U=(u_1,\cdots,u_N)\in (H_{0,r}^1(\Omega))^N$, the product space of $N$-copies of radially symmetric functions in $H^1_0$.
To avoid confusion of notations, we may also sometimes write $U= (U_1, \cdots, U_N)$.

In the following, for fixed $N$, we pick up a prime number $p$ and write $$N= Bp+R$$ for integers $B\geq 1$ and $R\geq 0$.
With this setting, the $N$ components of $U$ are divided into the union of $B$ subgroups with each group containing $p$ components and the remaining $R$ components:
\[
U=(U_1, ..., U_p, U_{p+1}, ..., U_{2p}, ......, U_{(B-1)p+1}, ..., B_{Bp}, U_{Bp+1},..., U_{Bp+R}).
\]
The solutions we construct exhibit qualitative properties for each component, between components within each groups, and between components from differen groups. Again for a continuous radial function $u$, $n(u)$ denotes its number of zeros whenever it makes sense.
Inductively we may define the following sequence of integers for $s=1,2,...$
\begin{equation}\label{def:Ks}
  \left\{
   \begin{array}{lr}
     K_{s+1}=8(p-1)^2\sum_{b=1}^B (P_b+1)\cdot K_s+5B(p-1)^2+(p-1)\sum_{b=1}^B (P_b +1)+1,\\
     K_1=8p(p+1)\sum_{b=1}^B (P_b +1)+5Bp^2.
   \end{array}
   \right.
   \end{equation}
Here is the first of our main results.
\begin{theorem}\label{t:existence}
Assume $\beta\leq-1$. Let $N=Bp+R$ with $p$ being a prime number, $B\geq 1$ and $R\geq 0$ integers. Fix non-negative integers $P_1,\cdots,P_B$ and $Q_1\cdots,Q_R$. Then there exist an infinite sequence of solutions $\{U_s= ((U_s)_1, \cdots, (U_s)_N)\}_{s=1}^\infty$ to Problem (\ref{e:AAAA}) such that
for all $s\geq 1$,
\begin{itemize}
  \item [$(1).$] (componentwise-prescribed nodes, 1). For any $b=1,\cdots,B$ and any $i=1,\cdots,p$, $$n((U_s)_{(b-1)p+i})=P_b;$$
  \item [$(2).$] (componentwise-prescribed nodes, 2). For any $r=1,\cdots,R$, $$n((U_s)_{Bp+r})=Q_r;$$
  \item [$(3).$] (Comparison within each group) for any $b=1,\cdots,B$ and any $j_1,j_2=1,\cdots,p$ with $j_1\neq j_2$,
      $$n((U_s)_{(b-1)p+j_1} - (U_s)_{(b-1)p+j_2})\in[ (P_b+1) \cdot K_{s}+2,(P_b+1) \cdot K_{s+1}+1].$$
      Here, $K_s$ is defined as in (\ref{def:Ks});
  \item [$(4).$] (Comparison between groups, 1) for any $b_1,b_2=1,\cdots,B$ and any $j_1,j_2=1\cdots,p$ with $b_1\neq b_2$, $$n((U_s)_{(b_1-1) p+j_1}-(U_s)_{(b_2 -1)p+j_2})\leq P_{b_1}+P_{b_2}+1;$$
  \item [$(5).$] (Comparison between groups, 2) for any $b=1,\cdots,B$, $j=1,\cdots,p$ and $r=1,\cdots,R$, $$n((U_s)_{(b-1)p+j}-(U_s)_{Bp+r})\leq P_b +Q_r +1;$$
  \item [$(6).$] (Comparison between groups, 3) for any $r_1 ,r_2=1,\cdots,R$ with $r_1\neq r_2$, $$n((U_s)_{Bp+r_1} - (U_s)_{Bp+r_2})\leq Q_{r_1}+Q_{r_2}+1.$$
\end{itemize}
\end{theorem}
\begin{remark}
1). We note the first two assertions in the above theorem give the prescribed nodal numbers for each component.\\
2). The last four assertions provide qualitative estimates of nodal data for differences between different components with some upper bounds and some lower bounds.\\
3).
In particular, the assertion (3) shows that this is an unbounded sequence of solutions since $K_s\to \infty$ as $s\to \infty$.\\
4).
The comparisons between different components involved in Theorem \ref{t:existence} are divided into four cases. The first case is the comparison in each group (Assertion (3)), i.e.,
  \begin{equation}
    \begin{tikzpicture}[>=stealth,baseline,anchor=base,inner sep=0pt]
      \matrix (foil) [matrix of math nodes,nodes={minimum height=0.5em}] {
        ( & u_1 & , & \cdots & , & u_p & ;\; & \cdots & ;\; & u_{(B-1)p+1} & , & \cdots & , & u_{Bp} & ; \;& u_{Bp+1} & , & \cdots & , & u_N ). \nonumber\\
      };
      \path[<->] ($(foil-1-2.north)+(0,1ex)$)   edge[blue,bend left=45]    ($(foil-1-6.north)+(0,1ex)$)             ;
    \end{tikzpicture}
  \end{equation}
The second case is the comparison between groups (Assertion (4)), i.e.,
  \begin{equation}
    \begin{tikzpicture}[>=stealth,baseline,anchor=base,inner sep=0pt]
      \matrix (foil) [matrix of math nodes,nodes={minimum height=0.5em}] {
        ( & u_1 & , & \cdots & , & u_p & ; & \cdots & ; & u_{(B-1)p+1} & , & \cdots & , & u_{Bp} & ; & u_{Bp+1} & , & \cdots & , & u_N ). \nonumber\\
      };
      \path[<->] ($(foil-1-2.north)+(0,1ex)$)   edge[blue,bend left=45]    ($(foil-1-10.north)+(0,1ex)$)             ;
    \end{tikzpicture}
  \end{equation}
The third case is the comparison between one component in groups and one of the rest of the components (Assertion (5)), i.e.,
  \begin{equation}
    \begin{tikzpicture}[>=stealth,baseline,anchor=base,inner sep=0pt]
      \matrix (foil) [matrix of math nodes,nodes={minimum height=0.5em}] {
        ( & u_1 & , & \cdots & , & u_p & ; & \cdots & ; & u_{(B-1)p+1} & , & \cdots & , & u_{Bp} & ; & u_{Bp+1} & , & \cdots & , & u_N ). \nonumber\\
      };
      \path[<->] ($(foil-1-2.north)+(0,1ex)$)   edge[blue,bend left=45]    ($(foil-1-16.north)+(0,1ex)$)             ;
    \end{tikzpicture}
  \end{equation}
  The forth case is the comparison between the rest of the components (Assertion (6)), i.e.,
  \begin{equation}
    \begin{tikzpicture}[>=stealth,baseline,anchor=base,inner sep=0pt]
      \matrix (foil) [matrix of math nodes,nodes={minimum height=0.5em}] {
        ( & u_1 & , & \cdots & , & u_p & ; & \cdots & ; & u_{(B-1)p+1} & , & \cdots & , & u_{Bp} & ; & u_{Bp+1} & , & \cdots & , & u_N & ). \nonumber\\
      };
      \path[<->] ($(foil-1-16.north)+(0,1ex)$)   edge[blue,bend left=45]    ($(foil-1-20.north)+(0,1ex)$)             ;
    \end{tikzpicture}
  \end{equation}
\end{remark}



Theorem \ref{t:existence}, while showing multiplicity, provide some upper bounds for the inter-componentwise comparisons. Conversely, the following result gives some lower bounds.

\begin{theorem}\label{t:nonexistence}
For Problem (\ref{e:AAAA}) with $\Omega$ radial, assume $\beta\leq-1$. Let $(u_1,\cdots,u_N)$ be a non-trivial radial solution to Problem (\ref{e:AAAA}). For any $i_1,i_2=1,\cdots,N$ with $i_1\neq i_2$, the following results hold.
\begin{itemize}
  \item [$(1)$]  $n(u_{i_1} - u_{i_2})\geq 1$ when $n(u_{i_1})=n(u_{i_2})=0$.
  \item [$(2)$] $n(u_{i_1} - u_{i_2})\geq\Big[\frac{\min\{n(u_{i_1}),n(u_{i_2})\}-1}{2}\Big]$
  when one of $n(u_{i_1})$ and $n(u_{i_2})$ is non-zero. Here, $[a]$ is the integer part of $a$.
\end{itemize}
\end{theorem}

\begin{remark}
To prove Theorem \ref{t:000}, we use two decompositions of $N=3$ and $P_1=\cdots=P_B=Q_1=\cdots=Q_R=0$.
Note that the integer $N=3$ can be decomposed in two cases.\\
(Case A): $p=2$, $B=1$, $R=1$.\\(Case B): $p=3$, $B=1$, $R=0$.\\
Applying Theorem \ref{t:existence} and Theorem \ref{t:nonexistence} with (Case B) we obtain the first three infinite sequence of positive solutions of Theorem \ref{t:000}, and
applying Theorem \ref{t:existence} and Theorem \ref{t:nonexistence} with (Case A) we obtain the fourth infinite sequence of positive solutions of Theorem \ref{t:000}.
\end{remark}


From a methodological point of view, the main tool in this paper is a variant of the symmetric mountain-pass theorem with proper invariant sets built in under the associated parabolic flow. The variational theory for the nodal critical points has been extensively studied, cf. \cite{CMT2000,LiuSun2001,Bartsch2001,BartschLiu2004,AckermanBartsch2005,LiuLiuWang2015} and the references therein.
To obtain the nodal critical points, the authors developed new frameworks by combining the classical mountain pass theorem and maximum principle.
The maximum principle guarantees the invariance of the positive cone under the negative gradient flow. This is shown in \cite{LiuSun2001,Bartsch2001,BartschLiu2004,LiuLiuWang2015}.
However, when considering the radial solutions with prescribed number of nodes, the parabolic flow is more useful here, cf. \cite{CMT2000,AckermanBartsch2005}. This is because of the well-known property of the parabolic equation in one spatial dimension or in the presence of radially symmetry that the number of nodes cannot increase along a flow line.

In our previous work \cite{LiWang2021}, we built a framework by embedding the parabolic flow into the $\mathbb{Z}_p$-symmetric mountain pass theorems. The fine nodal property (cf. Proposition \ref{p:L4} and Lemma \ref{l:NodalNumber} in this paper) are applied to prevent the nodal bumps from vanishing. In this paper, we further develop this framework by involving the comparisons of nodal data between the components of vector solutions and further reveal the structure of radial solutions to Problem (\ref{e:AAAA}).


\subsection{The idea and the organization of this paper}
This paper focuses on a symmetric mountain-pass theorem with certain invariant sets built in under the associated parabolic flow. Our paper is divided into two main parts: one part on the dynamics of parabolic equations and another part on the proof of the existence of solutions with certain properties by the dynamical contents. To be more precise, we outline our work as follows.
\begin{itemize}
  \item [$(1).$] In Section 2, we introduce the existence, regularity, global existence, the boundedness and nodal properties of the solutions to the corresponding parabolic problem. While giving proofs for some results here we refer most results to \cite{LiWang2021};
  \item [$(2).$] The proof of our main result Theorem \ref{t:existence} is given in Section 3. The proof is divided into three steps: $(a).$ For the set $\mathcal{A}_{\mathop{P}\limits^{\rightarrow},\bf{M}}^{\leq}\backslash\mathcal{D}$ (cf. (\ref{WorkingSet5}) and (\ref{Set:D})) of the vector-valued functions with prescribed upper bounds of comparison, we prove it admits a finite genus. See Lemma \ref{l:CompareSet}. $(b).$ We construct a sequence of sets $G_K\cap\partial\mathcal{A}\backslash(\bigcup_{j,q}C_{j,q}\cup H)$ (cf. (\ref{def:GK}), (\ref{def:A}), (\ref{def:Cjq1}-\ref{def:Cjq2}) and (\ref{def:H})) with componentwisely prescribed number of nodes with unbounded genus consisting of the vector-valued functions. See Corollary \ref{coro:lowerDk}. $(c).$ For a fixed set $\mathcal{A}_{\mathop{P}\limits^{\rightarrow},\bf{M}}^{\leq}\backslash\mathcal{D}$, we can always find a set $G_K\cap\partial\mathcal{A}\backslash(\bigcup_{j,q}C_{j,q}\cup H)$ with large genus which leads to an equilibrium point of the parabolic flow outside of $\mathcal{A}_{\mathop{P}\limits^{\rightarrow},\bf{M}}^{\leq}\backslash\mathcal{D}$. In this way, we find a solution to Problem (\ref{e:AAAA}) with prescribed number of nodes and desired comparison properties. To this end,
      \begin{itemize}
        \item [$(2.1).$] In Subsection \ref{Subsection:workingset}, we deal with set $\mathcal{A}_{\mathop{P}\limits^{\rightarrow},\bf{M}}^{\leq}\backslash\mathcal{D}$.  See Lemma \ref{l:CompareSet};
        \item [$(2.2).$] In Subsections \ref{Subsection:Simplex} and \ref{Subsection:nonvanishing}, we consider the set $G_K\cap\partial\mathcal{A}\backslash(\bigcup_{j,q}C_{j,q}\cup H)$. See Corollary \ref{coro:lowerDk};
        \item [$(2.3).$] In Subsections \ref{Subsection:Comparison} and \ref{Subsection:proof}, we find a  equilibrium point outside of $\mathcal{A}_{\mathop{P}\limits^{\rightarrow},\bf{M}}^{\leq}\backslash\mathcal{D}$. This is done by a comparison of genus;
      \end{itemize}
  \item [$(3).$] Theorem \ref{t:nonexistence} is proved in Section 4.
\end{itemize}


\vs

\noindent{\bf Notations.} Throughout this paper, generally, for a Banach space $X$, we write its norm as $\|\cdot\|_X$. Especially, for the Lebesgue space $L^p(\Omega)$ and Sobolev space $H^1_{0,r}(\Omega)$, the norms are denoted by $|\cdot|_p$ and $\|\cdot\|$, respectively. $\Omega$ is a bounded radial domain. To be precise, $\Omega$ is either a ball or an annulus.


\section{The parabolic settings}
Consider the following parabolic problem
\begin{equation}\label{e:BBBB}
    \left\{
   \begin{array}{lr}
     \partial_t u_j-{\Delta}u_{j}+ u_{j}= u^{3}_{j}+\beta\sum_{i=1, i\neq j}^N u_{j}u_{i}^{2} \,\,\,\,\,\,\,  \mbox{for}\ x\in \Omega\mbox{ and }t>0 ,\\
     u_{j}(0,x)=u_{j,0}(x)\in H_{0,r}^{1}(\Omega), \,\,\,\,\,\,\,\,j=1,\dots,N.
   \end{array}
   \right.
\end{equation}
Here, the constants $\beta\leq-1$. The domain $\Omega\subset\mathbb{R}^n$ is radial and $n=2,3$. In the following, we will use $\eta^t(U)$ or $(u_1(t), \cdots,u_N(t))$ to denote the solution to Problem (\ref{e:BBBB}). The solution is defined on the time interval $[0,T(U))$ with $T(U)$ denotes the maximal existing time of the solution. Especially, we will make use of the notion of $\omega$-set $\omega(A)$ of a subset $A$. To be precise, for a function $U\in (H_{0,r}^1(\Omega))^N$ with $T(U)=\infty$,
\begin{align}
\omega(U)=\{V\in (H_{0,r}^1(\Omega))^N|\exists t_n\to+\infty\mbox{ s.t. }\eta^{t_n}(U)\to V\mbox{ as }n\to\infty\}.\nonumber
\end{align}
For a subset $A\subset (H_{0,r}^1(\Omega))^N$,
\begin{align}
\omega(A)=\{V\in (H_{0,r}^1(\Omega))^N|\exists t_n\to+\infty\mbox{ and }\exists U_n\in A\mbox{ s.t. }\eta^{t_n}(U_n)\to V\mbox{ as }n\to\infty\}.\nonumber
\end{align}
We refer \cite{HenryBook} for more on the notions of dynamical systems.

As previously noted, the main tool in this paper is a mountain-pass theorem with the parabolic flow built-in.
To this end, we will investigate the related properties of Problem (\ref{e:BBBB}) in this section.

\subsection{Basic settings and results}
We provide a list of basic properties of Problem (\ref{e:BBBB}) which can be found in \cite{CMT2000,LiWang2021}. We only give the statements without proofs. Complete proofs for these results can be found in \cite{Chang2002,DM,HenryBook}.
We firstly address results on the local existence, uniqueness, regularity and continuous dependence.
\begin{theorem}\label{t:HeatExist}
For any initial value $U\in(H_{0,r}^{1}(\Omega))^N$, there is a unique solution $\eta^{t}(U)=(u_1(t),\cdots,u_n(t))$ to Problem (\ref{e:BBBB}) defined on its maximum interval $[0,T(U))$, satisfying
\begin{itemize}
  \item [$(I)$]$
   \eta^{t}(U)\in  C^1((0,T(U)),(L^2(\Omega))^N) \cap C([0,T(U)),(H_{0,r}^1(\Omega))^N)$;
  \item [$(II)$] for any $U\in (H_{0,r}^{1}(\Omega))^{2}$ and any $\delta\in[0,T(U))$, there are positive constants $r,K$ such that for any $t\in[0,\delta]$
  $$\|U-V\|_{(H_{0,r}^{1}(\Omega))^{N}}<r\,\,\,\,\Rightarrow\,\,\,\, \|\eta^{t}(U)-\eta^{t}(V)\|_{(H_{0,r}^1(\Omega))^N}\leq K\|U-V\|_{(H_{0,r}^{1}(\Omega))^{N}};$$
  \item [$(III)$] the trivial solution $(0,\cdots,0)\in (H_{0,r}^1(\Omega))^N$ is asymptotically stable in $(H_{0,r}^1(\Omega))^N$.
\end{itemize}
\end{theorem}

Furthermore, we can prove that $\eta^t(U)$ is an energy decreasing flow.
\begin{proposition}\label{p:L2partial_t}
For a solution $U(t)=(u_{1}(t),\dots,u_{N}(t))$ to Problem (\ref{e:BBBB}), we have
$$\frac{\partial}{\partial t}I(U(t))=-\sum_{j=1}^{N}\int|\partial_{t}u_{j}|^2\leq0.$$
\end{proposition}

\begin{corollary}
Problem (\ref{e:BBBB}) is dissipative.
\end{corollary}

Now we define the attracting domain of the trivial function in Soboelv space. The boundary of this domain will reveal a variety of topological structures.

\begin{corollary}\label{c:A}
Let
\begin{align}\label{def:A}
\mathcal{A}=\big\{U\in(H_{0,r}^{1}(\Omega))^{N}|T(U)=\infty\,\,\,\mbox{and}\,\, \lim_{t\to+\infty}\eta^{t}(U)=\theta\,\,\,\mbox{in}\,\,(H_{0,r}^1(\Omega))^N\big\}.
\end{align}
Then $\mathcal{A}$ is invariant under the heat flow and is open in $(H_{0,r}^{1}(\Omega))^{N}$. And both of $\mathcal{A}$ and $\partial\mathcal{A}$ are invariant under the flow generated by Problem (\ref{e:BBBB}). Moreover, $\inf_{U\in\partial\mathcal{A}}I>0$.
\end{corollary}

\subsection{Global behaviours}
In this subsection, we focus on the global existence and boundedness on $\partial\mathcal{A}$.
First, we establish the global existence of orbits on $\partial\mathcal{A}$ using the method in Cazenave and Lions \cite{CazenaveLions1984,Chang2002}.
\begin{lemma}\label{l:globalexists}
Suppose $\left\{I(\eta^t(U))\right\}_{t\in[0,T(U))}$ is bounded from below. Then $U(t)$ exists globally in $(H_{r,0}^1(\Omega))^N$.
\end{lemma}

\begin{corollary}
For any $U\in\partial\mathcal{A}$, $T(U)=\infty$.
\end{corollary}

And if we improve the regularity of the initial data, we can obtain a global $H^1$-boundedness on $\partial\mathcal{A}$.

\begin{proposition}\label{prop:GlobalBound}
For any $U\in\partial\mathcal{A}\cap (H^2_{0,r}(\Omega))^N$, $\|\eta^{t}(U)\|_{(H_{0,r}^1(\Omega))^N}\leq C$ for any $t\geq0$. Here, the constant $C>0$ depends continuously on the $H^2$-norm of the  initial data.
\end{proposition}

By the variation of constant, we can establish the $H^2$-boundedness of the orbit.

\begin{corollary}\label{coro:LinftyBdd}
For any $U\in\overline{\mathcal{A}}\cap(H^2_{r}(\Omega))^N$, there is a constant $C(U)>0$ such that for any $t\geq0$, $|\eta^t(U)|_\infty\leq C(U)$. Furthermore, the constant $C(U)$ can be taken as a continuous function of $\|U\|_{(H_{r}^2(\Omega))^N}$.
\end{corollary}
\noindent{\bf Proof.}
Using the variation of constant, we get
\begin{align}
\|U(t)\|_{(H^{\frac{4}{3}}(\Omega))^N}&\leq Ce^{- t}\|U\|_{(H^{\frac{4}{3}}(\Omega))^N}+C\int_{0}^{t} \frac{e^{-(t-s)}}{|t-s|^{\frac{2}{3}}}\Big|\sum_{i=1}^N\sum_{j=1}^N\beta_{ij}u_i u_j^2\Big|_{2}ds\nonumber\\
&\leq Ce^{- t}\|U\|_{(H^{\frac{4}{3}}(\Omega))^N}+ C\cdot\sup_{t\geq0}|U(t)|_{6}^3\int_{0}^t\frac{e^{-(t-s)}}{|t-s|^{\frac{2}{3}}}ds\nonumber\\
&\leq Ce^{-t}\|U\|_{(H^{\frac{4}{3}}(\Omega))^N}+ C\cdot\sup_{t\geq0}\|U(t)\|_{(H^1(\Omega))^N}^3\int_{0}^t\frac{e^{-(t-s)}}{|t-s|^{\frac{4}{3}}}ds\nonumber\\
&\leq Ce^{- t}\|U\|_{(H^{\frac{4}{3}}(\Omega))^N}+ C\cdot\sup_{t\geq0}\|U(t)\|_{(H^1(\Omega))^N}^3\int_{0}^\infty\frac{e^{- s}}{|s|^{\frac{2}{3}}}ds\leq C_1.\nonumber
\end{align}
Therefore, for any $U\in\overline{\mathcal{A}}\cap(H^2(\Omega))^N$, there is a constant $C_1>0$ such that for any $t\geq0$, $\|U(t)\|_{(H^{\frac{4}{3}}(\Omega))^N}\leq C_1$. We refer \cite[Remark 2.2.10]{LuBook} for the estimate on the heat kernel between fractional spaces. With a similar approach,
\begin{align}
\|U(t)\|_{(H^{\frac{5}{3}}(\Omega))^N}&\leq Ce^{- t}\|U\|_{(H^{\frac{5}{3}}(\Omega))^N}+ C\int_{0}^{t}\frac{e^{-(t-s)}}{|t-s|^{\frac{2}{3}}}\Big\|\sum_{i=1}^N\sum_{j=1}^N\beta_{ij}u_i u_j^2\Big\|_{(H^{\frac{1}{3}}(\Omega))^N}ds\nonumber\\
&\leq Ce^{-t}\|U\|_{(H^{\frac{5}{3}}(\Omega))^N} +C\cdot\sup_{t\geq0}\|U(t)\|_{(W^{\frac{1}{3},6}(\Omega))^N}^3 \int_{0}^t\frac{e^{-(t-s)}}{|t-s|^{\frac{2}{3}}}ds.\nonumber
\end{align}
\begin{flushright}
$\Box$
\end{flushright}

\begin{remark}
1). We choose a $\frac{1}{3}$-bootstrap in order to avoid Lions-Magenes space (cf. \cite[Chapter 33]{TartarBook}).

2). Furthermore, we can prove the following general boundedness result.
For any $U\in\overline{\mathcal{A}}\cap(H^k(\Omega))^N$ and for any integer $k$, there is a constant $C(U, k)>0$ such that for any $t\geq0$, $\|\eta^t(U)\|_{H^k(\Omega)}\leq C(U,k)$. Furthermore, the constant $C(U,k)$ can be taken as a continuous function of $\|U\|_{(H^k(\Omega))^N}$.
The proof is similar to the one of Corollary \ref{coro:LinftyBdd} and we omit it here.
\end{remark}

\subsection{Nodal properties}

In this paper, we need two kind of nodal properties. They each describe
\begin{itemize}
  \item [$(1).$] The number of bumps;
  \item [$(2).$] Then size of bumps.
\end{itemize}

To begin with, it is necessary to introduce the notions of nodal number and bump.
\begin{definition}
For a continuous radial function $u:\Omega\to\mathbb{R}$, we define the number of nodes of the function $u$ to be the the largest number $k$ such that there exist a sequence of real numbers $x_0,\cdots,x_k$ such that $0<x_{0}<x_{1}<\dots<x_{k}$ and
      $$u|_{|x|=x_j}\cdot u|_{|x|=x_{j+1}}<0,\,\,\,\,\,\,j=0,\dots,k-1.$$
Denote the nodal number of the function $u$ by $n(u)$. We define its $q$-th bump $q=1,...,k+1$ by
      \begin{align}
      u_{1}(x)&=\chi_{\{\mbox{sgn}(u(x))=\mbox{sgn}u(x_0)\}}\cdot\chi_{\{|x|<x_{1}\}}\cdot u(x),  \nonumber\\
      u_{q}(x)&=\chi_{\{\mbox{sgn}(u(x))=\mbox{sgn}u(x_{q-1})\}}\cdot\chi_{\{x_{q-2}<|x|<x_{q}\}}\cdot u(x),\,\,\,\,q=2,\dots,k,\nonumber\\
      u_{k+1}(x)&=\chi_{\{\mbox{sgn}(u(x))=\mbox{sgn}u(x_q)\}}\cdot\chi_{\{x_{q-1}<|x|\}}\cdot u(x).\nonumber
      \end{align}
      For the $j$-th component $u_{j}$ of the vector-valued function $U=(u_{1},\dots,u_{N})$, we denote its $q$-th bump by $u_{j,q}$.
\end{definition}

Now we present a result from \cite{IshiwataLi,WeiWeth2008}, which concerns the number of bumps.
\begin{lemma}\label{l:NodalNumber}
For the problem
\begin{equation}\label{e:LINEAR}
    \left\{
   \begin{array}{lr}
     \partial_t w-{\Delta}w+g(x,t)w=0 \,\,\,\,\,\,\,\,  \mbox{in}\ \Omega ,\\
     w(0,x)=w_0(x)
   \end{array}
   \right.
\end{equation}
with radial initial data and $g(x,t)>0$ is radial in $x$, the nodal number $n(w(x,t))$ of the classical solution $w(x,t)\in H^1_{r}$ is non-increasing as the time $t\geq0$ increases. 
\end{lemma}
This is a classical result also can be found in \cite{Matano1982} and the references therein.
For any constant $\Lambda_j$ and $v_j(t,x)=e^{-\Lambda_j t}$ satisfies
$$\partial_t v_j-{\Delta}v_j+(\Lambda_j+1- u_j^3-\beta\sum_{i\neq j}u_i^2)v_j=0 \mbox{ in }\ \Omega.$$
For large $\Lambda_j>0$, $\Lambda_j+1- u_j^3-\beta\sum_{i\neq j}u_i^2>0$ for any $(x,t)\in\Omega\times[t_0,t_1]$. Then, the nodal number of $v_j$ is non-increasing in $[t_0,t_1]$. This means that we prove the following corollary.
\begin{corollary}\label{c:nodalnumber}
For a solution $(u_1(t),\cdots,u_N(t))$ to Problem (\ref{e:BBBB}), for any $j=1,\cdots,N$, $n(u_j(t))$ is non-increasing.
\end{corollary}

Then,
we consider a property presented by \cite{CMT2000,LiWang2021}. This property claims that the parabolic flow preserves the smallness of the $L^4$-norm of the bumps.
\begin{proposition}\label{p:L4}
Assume that in Problem (\ref{e:BBBB}) $\beta\leq-1$.
There is a positive number $\rho>0$ such that for any $U(t)=(u_1(t),\cdots,u_N(t))$ solving Problem (\ref{e:BBBB}), if $|u_{j,q}(0)|_{4}<\rho$ then $|u_{j,q}(t)|_{4}<\rho$ for $t\geq0$ and $|u_{j,q}(t)|_{4}>0$ for any $j=1,\cdots,N$ and $q=1,\cdots,n(u_j)+1$.
\end{proposition}
For the sake of completeness, we present the proof of Proposition \ref{p:L4} here. To do this, we need the following result on the differentiation on the time scale.
\begin{lemma}\label{l:diff}
For a solution $(u_1,\cdots,u_N)$ to Problem (\ref{e:BBBB}) and any $t_0\in(0,T(U))$, a bump $u_{j,p}(t_0)$ satisfies that
$\partial_t\int|u_{j,p}|^4\Big|_{t=t_0}=4\int\partial_t u_j\cdot (u_{j,p})^3\Big|_{t=t_0}$.
\end{lemma}
\noindent{\bf Proof.}
Without loss of generality, let us assume that there are two nodes $a(t)<b(t)$ with $t\in[t_0,t_0+\Delta t]$
such that
\begin{itemize}
  \item $\partial_t\int|u_{j,p}|^4\Big|_{t=t_0}=\partial_t\int_{a(t_0)}^{b(t_0)}|u|^4$;
  \item $a(t_0+\Delta t)\to a(t_0)$ and $b(t_0+\Delta t)\to b(t_0)$ as $\Delta t\to 0^+$.
\end{itemize}
Let us begin by definition of derivatives.
\begin{align}
&\frac{1}{\Delta t}  \Big[\int_{a(t_0+\Delta t)}^{b(t_0+\Delta t)}|u(t_0+\Delta t,x)|^4dx - \int_{a(t_0)}^{b(t_0)}|u(t_0,x)|^4dx\Big]\nonumber\\
=&\frac{1}{\Delta t} \Big[\int_{a(t_0+\Delta t)}^{b(t_0+\Delta t)} - \int_{a(t_0)}^{b(t_0)}|u(t_0+\Delta t,x)|^4dx\Big]\nonumber \\
& +\frac{1}{\Delta t} \Big[\int_{a(t_0)}^{b(t_0)}|u(t_0+\Delta t,x)|^4-|u(t_0,x)|^4dx\Big]\nonumber\\
=&:I_1+I_2.\nonumber
\end{align}
Here,
\begin{align}
I_1=&\frac{1}{\Delta t} \Big[\int_{a(t_0+\Delta t)}^{b(t_0+\Delta t)} - \int_{a(t_0)}^{b(t_0)}|u(t_0+\Delta t,x)|^4dx\Big]\nonumber\\
=&\;\;\;\;\; \frac{1}{\Delta t} \Big[\int_{b(t_0)}^{b(t_0+\Delta t)} - \int_{a(t_0)}^{a(t_0+\Delta t)}|u(t_0,x)|^4dx\Big] + \frac{1}{\Delta t} \Big[\int_{a(t_0+\Delta t)}^{b(t_0+\Delta t)} - \int_{a(t_0)}^{b(t_0)}o(1)dx\Big]\nonumber\nonumber\\
=&|u(t_0,x)|^4\Big|_{x=b(t_0)}-|u(t_0,x)|^4\Big|_{x=a(t_0)}+o(1)\to0\nonumber
\end{align}
as $\Delta t\to0$. The second line is due to the continuity of $t\mapsto u(x,t)$ and the last line is due to the differentiation of the integral, cf. \cite[Corollary 1/p.5]{SteinBook}.
And,
\begin{align}
I_2\to\partial_t\int_{a(t_0)}^{b(t_0)}|u|^4=\partial_t\int|u_{j,p}|^4\nonumber
\end{align}
holds obviously.

\begin{flushright}
$\Box$
\end{flushright}

Based on this, we have

\vs

\noindent{\bf Proof of Proposition \ref{p:L4}.}
According to the assumptions in Proposition \ref{p:L4}, there exists a small $\varepsilon>0$ such that if
$$(-1)^{q+1}u_{j}(x_{q},0)>0,$$
then
$$(-1)^{q+1}u_{j}(x_{q},t)>0$$
for any $t\in[0,\varepsilon]$. Hence, due to the definition of bump $u_{j,q}$, the differential $\frac{\partial}{\partial t}\int|u_{j,q}|^{4}$ is well-defined. By Lemma \ref{l:diff},
\begin{align}
\frac{\partial}{\partial t}\int|u_{j,q}|^{4} &=4\int u_{j,q}^{3}\partial_{t}u_{j,q}=4\int u_{j,q}^{3}\partial_{t}u_{j}\nonumber\\
&=4\int u_{j,q}^{3}\Big(\Delta u_{j}- u_{j}+ u_{j}^{3}+\sum_{i\neq j}\beta u_{j}u_{i}^{2}\Big)\nonumber\\
&=-3\int |\nabla(u_{j,q}^{2})|^{2}-4 \int u_{j,q}^{4}+4 \int u_{j,q}^{6}+4\sum_{i\neq j}\beta \int u_{j,q}^{4}u_{i}^{2}.\nonumber
\end{align}
Denote $W=u_{j,q}^{2}$. Noticing
$$\frac{1}{3}=\frac{\frac{1}{2}}{6}+\frac{1-\frac{1}{2}}{2},$$
we have from Sobolev embedding,
$$|W|_{3}^{3}\leq C\|W\|^{\frac{3}{2}}|W|_{2}^{\frac{3}{2}}.$$
Therefore,
\begin{align}
\frac{\partial}{\partial t}\int |u_{j,p}|^{4}&\leq-C\|W\|^{2}+C\|W\|^{\frac{3}{2}}|W|^{\frac{3}{2}}_2\nonumber\\
&\leq-C\|W\|^{\frac{3}{2}}|W|_{2}^{\frac{1}{2}}+C\|W\|^{\frac{3}{2}}|W|^{\frac{3}{2}}_2\nonumber\\
&=-C\|W\|^{\frac{3}{2}}|W|_{2}^{\frac{1}{2}}\big(1-C|W|_{2}\big)\nonumber\\
&=-C\|W\|^{\frac{3}{2}}|W|_{2}^{\frac{1}{2}}\big(1-C|u_{j,q}|_{4}^{2}\big)<0\nonumber
\end{align}
for $|u_{j,q}|_{4}$ small enough.
\begin{flushright}
$\Box$
\end{flushright}

\subsection{A Unique Continuation Result for A Backward Parabolic Inequality}
To compute the linking structure, we need one more basic property of the parabolic flow. This version of unique continuation for a backward parabolic inequality is a special case of \cite[Theorem 1.(i)]{Fernandez2003}.
To state the next theorem, we need the following two function spaces:
\begin{itemize}
  \item $L_t^{1}\left((0,T),L^{\infty}(\Omega)\right)=\Big\{u:(0,T)\times\Omega\to\mathbb{R}\Big|
      \int_0^T|u|_{L^q(\Omega)}(t)dt<+\infty\Big\}$;
  \item $C_0(\Omega)=\{u\in C(\Omega)|\mbox{supp}u\mbox{ is compact in }\Omega\}$.
\end{itemize}

\begin{theorem}\label{t:backward}
Consider a function $V(t,x):[0,T]\times\Omega\to\mathbb{R}$ with $\|V\|_{L_t^{1}\left((0,T),L^{\infty}(\Omega)\right)}<\infty$. If a function $u$ satisfies the inequality $|\Delta u+\partial_t u|\leq V(t,x)|u|$ with $u(0,x)\in C_0(\Omega)$, then $u\equiv0$ in $[0,T]\times\Omega$.
\end{theorem}

\section{Proof of Theorem \ref{t:existence}}

\subsection{Outline of the proof of Theorem \ref{t:existence}}
Since our proof is rather long, we first sketch the idea in this subsection. Like in \cite{LiWang2021}, the main idea is to construct a variant of symmetric mountain-pass theorem on $\partial\mathcal{A}$.
To find a solution to Problem (\ref{e:AAAA}) is to find a equilibrium point of (\ref{e:BBBB}). To do this, we find a initial data $U\in\partial\mathcal{A}$ whose omega set contains a solution to Problem (\ref{e:AAAA}). Our argument in this section is a refinement of the above idea. To be precise,
\begin{itemize}
  \item [$(1).$] In Subsection \ref{Subsection:workingset}, we introduce some basic notions and several properties of certain sets;
  \item [$(2).$] For the set $\mathcal{A}_{\mathop{P}\limits^{\rightarrow},\bf{M}}^{\leq}\backslash\mathcal{D}$ (cf. (\ref{WorkingSet5}) and (\ref{Set:D})) of the vector-valued functions with prescribed upper bounds of comparison, we prove it admits a finite genus, cf. Lemma \ref{l:CompareSet} of Subsection \ref{Subsection:workingset}. Here, the matrix $\bf{M}$ contains the comparison among the groups;
  \item [$(3).$] We construct a sequence of sets consisting of the vector-valued functions with componentwisely prescribed number of nodes and each of them are non-vanishing under the parabolic flow. Such a set is proved to have a infinite genus, cf. Subsections \ref{Subsection:Simplex}, \ref{Subsection:nonvanishing} and \ref{Subsection:Comparison}. We prove this by proving it contains a sequence of sets $G_K\cap\partial\mathcal{A}\backslash(\bigcup_{j,q}C_{j,q}\cup H)$ (cf. (\ref{def:GK}), (\ref{def:A}), (\ref{def:Cjq1}-\ref{def:Cjq2}) and (\ref{def:H})) with unbounded genus. which is ensured by Corollary \ref{coro:lowerDk};
  \item [$(4).$] For the set  $\mathcal{A}_{\mathop{P}\limits^{\rightarrow},\bf{M}}^{\leq}\backslash\mathcal{D}$ with fixed $\bf{M}$, we can choose the set $G_K\cap\partial\mathcal{A}\backslash(\bigcup_{j,q}C_{j,q}\cup H)$ with large genus and an initial data in $G_K\cap\partial\mathcal{A}\backslash(\bigcup_{j,q}C_{j,q}\cup H)$ outside of $\mathcal{A}_{\mathop{P}\limits^{\rightarrow},\bf{M}}^{\leq}\backslash\mathcal{D}$ leading to an equilibrium point with comparison more than we prescribed in $\bf{M}$. See Subsection \ref{Subsection:proof}.
\end{itemize}

\vs

Due to Step $(3)$, the solution and its comparison among the components admit prescribed numbers of nodes.
Therefore, at the end, we only need to figure the relation with $\bf{R}$ and the constant $K$ in $G_K\cap\partial\mathcal{A}\backslash(\bigcup_{j,q}C_{j,q}\cup H)$. This is done in Subsection \ref{Subsection:proof}.

\subsection{Basic working spaces}\label{Subsection:workingset}
Let $P_1,\cdots,P_B,Q_1,\cdots,Q_R$ be the integers given in Theorem \ref{t:existence}. Here, $N=pB+R$ with $B>0$, $R\geq0$ and $p$ prime as we assumed. Denote $\mathop{P}\limits^{\rightarrow}:=(P_1,\cdots,P_B,Q_1,\cdots,Q_R)$.
The space of vector-valued functions with component-wisely prescribed number of nodes is denoted as
\begin{align}\label{WorkingSet1}
\mathcal{A}_{\mathop{P}\limits^{\rightarrow}}:=\Big\{U\in\overline{\mathcal{A}}\cap(H_{0,r}^2(\Omega))^N \Big|n(u_{(b-1)p+j})=P_b, & \,\,b=1,\cdots,B; \,j=1\cdots,p\nonumber\\
&\mbox{ and }n(u_{Bp+r})=Q_i,\,\,r=1,\cdots,R\Big\}.
\end{align}
In order to compare the nodal numbers of different components, we analyze all the possible differences within our setting. Recall that we apply a partial permutations in groups of components defined as
\begin{align}\label{PermutationPP}
\sigma(u_{1},u_{2},\dots,u_{p};&\dots\dots;u_{(B-1)p+1},u_{(B-1)p+2},\dots, u_{Bp};u_{Bp+1},\cdots,u_N)\nonumber\\
&=(u_{2},\dots,u_{p},u_{1};\dots\dots;u_{(B-1)p+2}, \dots,u_{Bp},u_{(B-1)p+1};u_{Bp+1},\cdots,u_N).
\end{align}
We first consider the comparisons among $u_{pb+1-p},\cdots,u_{pb}$ for each $b=1,\cdots,B$. Due to the nature of the permutation action, for each group, such as Group $b$, we have the following list of the possible couples:
\begin{itemize}
  \item [$\mbox{Case }(b,1).$]
  {\bf Distance=$1$ or $p-1$.} $(u_{(b-1)p+1},u_{(b-1)p+2})$, $(u_{(b-1)p+2}, u_{(b-1)p+3})$, $\cdots$, $(u_{bp-2},u_{bp-1})$, $(u_{bp-1},u_{bp})$, $(u_{pb},u_{pb+1-p})$;
  \item [$\mbox{Case }(b,2).$] {\bf Distance=$2$ or $p-2$.} $(u_{(b-1)p+1},u_{(b-1)p+3})$, $(u_{(b-1)p+2},u_{(b-1)p+4})$, $\cdots$, $(u_{bp-2},u_{bp})$, $(u_{bp-1},u_{(b-1)p+1})$, $(u_{bp},u_{(b-1)p+2})$;
  \item [$\cdots.$] {\bf $\cdots$.} $\cdots$.
  \item [$\mbox{Case }\left(b,\frac{p-1}{2}\right).$] {\bf Distance=$\frac{p-1}{2}$ or $\frac{p+1}{2}$.} $(u_{(b-1)p+1},u_{(b-1)p+\frac{p+1}{2}})$, $(u_{(b-1)p+2},u_{(b-1)p+\frac{p+3}{2}})$, $\cdots$, $(u_{bp-2},u_{(b-1)p+\frac{p-1}{2}-2})$, $(u_{bp-1},u_{(b-1)p+\frac{p-1}{2}-1})$, $(u_{bp},u_{(b-1)p+\frac{p-1}{2}})$.
\end{itemize}
Observe that there are $p$ couples for each case and in total there are $p\cdot\frac{p-1}{2}=C_p^2$ cases. It is worth to be pointed out that for the $p=2$ case only the $(u_1,u_2)$ case will occur. To proceed on with our discussion, we define the following set for Case $(b,q)$ with $b=1,\cdots,B$ and $q=1,\cdots,p$:
\begin{align}\label{WorkingSet2}
\mathcal{B}_{\mathop{P}\limits^{\rightarrow};b,q,M}^{=}:=\big\{U\in \mathcal{A}_{\mathop{P}\limits^{\rightarrow}}\cap\partial\mathcal{A} \Big|\mbox{for any couple }(u_i,u_j)\mbox{ from Case }(b,q)\mbox{ we have }n(u_i - u_j)=M\big\}
\end{align}
and
\begin{align}\label{WorkingSet3}
\mathcal{B}_{\mathop{P}\limits^{\rightarrow};b,q,M}^{\leq}:=\big\{U\in \mathcal{A}_{\mathop{P}\limits^{\rightarrow}}\cap\partial\mathcal{A}  \Big|\mbox{for any couple }(u_i,u_j)\mbox{ from Case }(b,q)\mbox{ we have }n(u_i - u_j)\leq M\big\}
\end{align}
For any positive integers $M_q^{(b)}:\,M_1^{(1)},\cdots, M_1^{(B)},\cdots, M_2^{(1)}, M_{\frac{p-1}{2}}^{(1)},\cdots ,M_{\frac{p-1}{2}}^{(B)}$ (only $M_1^{(1)},\cdots,M_{1}^{(B)}$ when $p=2$ or $3$), denote
\begin{align}\label{WorkingSet4}
\mathcal{A}_{\mathop{P}\limits^{\rightarrow},\bf{M}}^=:=\cap_{b=1}^B\cap_{q=1}^{\frac{p-1}{2}} \mathcal{B}_{\mathop{P}\limits^{\rightarrow};b,q,M_q^{(b)}}^{=}
\end{align}
and
\begin{align}\label{WorkingSet5}
\mathcal{A}_{\mathop{P}\limits^{\rightarrow},\bf{M}}^{\leq}:=\cap_{b=1}^B\cap_{q=1}^{\frac{p-1}{2}} \mathcal{B}_{\mathop{P}\limits^{\rightarrow};b,q,M_q^{(b)}}^{\leq}
\end{align}
Here,
\begin{align}\label{Matrix}
{\bf M}=\left(
  \begin{array}{cccc}
    M_1^{(1)} & M_2^{(1)} & \cdots & M_{\frac{p-1}{2}}^{(1)} \\
    M_1^{(2)} & M_2^{(2)} & \cdots & M_{\frac{p-1}{2}}^{(2)} \\
    \vdots & \vdots & \ddots & \vdots \\
    M_1^{(B)} & M_2^{(B)} & \cdots & M_{\frac{p-1}{2}}^{(B)} \\
  \end{array}
\right).
\end{align}
In this matrix, the $(b,q)$ element $M_q^{(b)}$ denotes the nodal number or the maximum of the nodal numbers of the difference of the couples in Case $(b,q)$.

Now we use the reduction in \cite{WeiWeth2008} to calculate the genus $\gamma_{pp}\Big(\mathcal{A}_{\mathop{P}\limits^{\rightarrow},\bf{M}}^{\leq}\Big)$. Here, we recall the definition of $\gamma_{pp}$, with $pp$ indicating partial-permutation.
Define
\begin{itemize}
  \item $\mbox{Fix}_{pp}=\{U\in(H_{0,r}^1(\Omega))^N|\sigma(U)=U\}$;
  \item $\mathcal{E}_{pp}=\{A\subset(H_{0,r}^1(\Omega))^N\backslash\mbox{Fix}_{pp}|\sigma(A)=A\mbox{ and }A\mbox{ is compact}\}$.
\end{itemize}
For any $A\in\mathcal{E}_{pp}$, we define the set
\begin{align}
\mathbb{I}_1(A):=\{m\in\mathbb{N}_+|\exists\mbox{ a continous }h:A\to\mathbb{C}^m\backslash\{0\}\mbox{ s.t. }h(\sigma U)=e^{\frac{2\pi i}{p}}h(U)\mbox{ }\forall U\in A\}.\nonumber
\end{align}
Here, the mapping $\sigma$ is defined as in (\ref{PermutationPP}).
Then, the index $\gamma_{pp}$ is defined as
\begin{equation}\label{def:GammaPP}
\gamma_{pp}(A)=\left\{
\begin{aligned}
&\min\{m|m\in\mathbb{I}_1(A)\}&\mbox{ if }\mathbb{I}_1(A)\neq\emptyset, \\
&\infty & \mbox{ if }\mathbb{I}_1(A)=\emptyset.
\end{aligned}
\right.
\end{equation}

Then, for any $A,B\in\mathcal{E}_{pp}$,
we have
\begin{proposition}\label{prop:genus}
\begin{itemize}
        \item [$(1).$] If $A\subset B$, then $\gamma_{pp}(A)\leq\gamma_{pp}(B)$;
        \item [$(2).$] $\gamma_{pp}(A\cup B)\leq\gamma_{pp}(A)+\gamma_{pp}(B)$;
        \item [$(3).$] if $g:A\to (H_{0,r}^2(\Omega))^N\backslash \mbox{Fix}_{\sigma}$ is continuous and satisfies
        $g(\sigma(u))=\sigma g(u)$ for all $ u\in A$,
        then $$\gamma_{pp}(A)\leq\gamma_{pp}(\overline{g(A)});$$
        \item [$(4).$] if $\gamma_{pp}(A)>1$, then $A$ is an infinite set;
        \item [$(5).$] if $A$ is compact and $\gamma_{pp}(A)<\infty$, then there exist an open $\sigma$-invariant neighbourhood $\mathcal{N}$ of $A$
        such that $\gamma_{pp}(A)=\gamma_{pp}(\overline{\mathcal{N}})$;
        \item [$(6).$] if $S$ is the boundary of a bounded neighbourhood of the origin in a $m$-dimensional complex linear space such that
       $e^{\frac{2\pi i}{p}}U \in S $ for any $U\in S$, and $\Psi:S\to (H_{0,r}^1(\Omega))^N\backslash \mbox{Fix}_\sigma $ is continuous and satisfies for any $U\in S$,
       $\Psi(e^{\frac{2\pi i}{p}}U)=\sigma(\Psi(U))$, then $\gamma_{pp}(\Psi(S))\geq m$;
       \item [$(7).$] Let $A$ be a closed set such that $A\subset (H_{0,r}^1(\Omega))^N\setminus \mbox{Fix}_\sigma$, and $\cap_{i=0}^{p-1} \sigma^i(A) =\emptyset $.
Then $\gamma_{pp}(\mathbb Z_p(A))\leq p-1$.
   \end{itemize}
\end{proposition}
All of the properties are standard but the last one. Readers can find Assertions (1-6) in \cite{TianWang2011}. \cite{RabinBook,Wang1989} are referred as general introductions.
We refer \cite{LiWang2021} for Assertion (7).

Besides these standard properties, there is an additional lemma to consider. This property applies to general $\mathbb{Z}_p$-genus.

\begin{lemma}\label{l:2genus}
For $i=1,2$, $E_i$ denote the Banach spaces equipped with $\mathbb{Z}_p$-actions $\sigma_i$. Denote $\mathcal{E}_i:= E_i\backslash\mbox{Fix}_{\sigma_i}$ with $\mbox{Fix}_{\sigma_i}$ denoting the sets of fixed points of $\sigma_i$, respectively. Suppose that for sets $A_i\subset\mathcal{E}_i$ with $\sigma_i(A_i)=A_i$, if there exists a map $f:A_1\to A_2$ with $f(\sigma_1 x)=\sigma_2 \circ f(x)$ for any $x\in A_1$, then we have $\gamma_{\sigma_1}(A_1)\leq\gamma_{\sigma_2}(A_2)$. Here, $\gamma_{\sigma_i}$ is the $\mathbb{Z}_p$-genus generated by the action $\sigma_i$ for $i=1,2$.
\end{lemma}


This is a direct consequence of result of Assertion (6).

Let us denote
\begin{align}\label{Set:D}
\mathcal{D}=\{U\in\partial\mathcal{A}\cap(H_{0,r}^2(\Omega))^N|\exists T>0\mbox{ s.t. }|\eta^T(U)_i|_4<\rho\mbox{ for some }i=1,\cdots,N\}.
\end{align}
Here, $\eta^T(U)_i$ is the $i$-th component of the vector-valued function $\eta^T(U)$ and $\rho$ is the constant in Proposition \ref{p:L4}.

\begin{lemma}\label{l:CompareSet}
It holds that
$$\gamma_{pp}\Big(\mathcal{A}_{\mathop{P}\limits^{\rightarrow},\bf{M}}^{\leq}\backslash\mathcal{D}\Big)\leq (p-1)\left(2\min_{b=1,..,B;q=1,...,\frac{p-1}{2}}M_q^{(b)}+3\right).$$
Here, the index $\gamma_{pp}$, the set $\mathcal{A}_{\mathop{P}\limits^{\rightarrow},\bf{M}}^{\leq}$ and the matrix $\bf{M}$ are defined in (\ref{def:GammaPP}), (\ref{WorkingSet5}) and (\ref{Matrix}).
\end{lemma}
\noindent{\bf Proof.}
We prove this result by a reduction.
Without loss of generality, let us assume that
\begin{align}
M_1^{(1)}=\min\Big\{M_q^{(b)}\Big|b=1,..,B;q=1,...,\frac{p-1}{2}\Big\}.\nonumber
\end{align}
And we will prove the lemma by reduction.
As a first step, we estimate the genus of $\gamma_{pp}\Big(\mathcal{A}_{\mathop{P}\limits^{\rightarrow},\bf{M}_0}^{\leq}\Big)$. Here, ${\bf M}_0$ is the matrix ${\bf M}$ with the $(1,1)$-element replaced by 0.
It is easy to see that
$$\mathcal{A}_{\mathop{P}\limits^{\rightarrow},\bf{M}_0}^\leq  \subset\bigcup_{i=0}^{p-1}\sigma^i \Big\{U\in\mathcal{A}_{\mathop{P}\limits^{\rightarrow}}\Big|u_1\geq u_2\Big\}.$$
Denote $A=\Big\{U\in\mathcal{A}_{\mathop{P}\limits^{\rightarrow}}\Big|u_1\geq u_2\Big\}$.
One can easily verify that $\cap_{i=0}^{p-1}\sigma^i(A)\subset \mathcal{D}$. Therefore, according to Assertion (7) of \ref{prop:genus}, we get $$\gamma_{pp}\Big(\mathcal{A}_{\mathop{P}\limits^{\rightarrow},\bf{M}_0}^{\leq}\Big)\leq p-1.$$

Now let us reduce on the $(1,1)$-element of the matrix ${\bf M}$.
Define
\begin{equation}
e_{ij}=\left\{
\begin{aligned}
1 & \mbox{  if} & (i,j)=(1,1), \nonumber\\
0 & \mbox{  if} & \mbox{otherwise}.\nonumber
\end{aligned}
\right.
\end{equation}
We will compare the difference between $\mathcal{A}_{\mathop{P}\limits^{\rightarrow},\bf{M}}^\leq$ and $\mathcal{A}_{\mathop{P}\limits^{\rightarrow},\bf{M}+{\bf e}}^\leq$. Observe that
\begin{align}
\mathcal{A}_{\mathop{P}\limits^{\rightarrow},\bf{M}+{\bf e}}^\leq=\mathcal{A}_{\mathop{P}\limits^{\rightarrow},\bf{M}}^\leq\cup \mathcal{A}_{\mathop{P}\limits^{\rightarrow},\bf{M}+{\bf e}}^= \cup A,
\end{align}
for some set $A$. Here, any $U\in A$, for the first group of the components $u_{1},\cdots,u_{p}$, we can find two couples of Case $(1,1)$, say $(i,j)$ and $(i',j')$ such that
\begin{itemize}
  \item $(i,j)\neq(i',j')$;
  \item $n(u_i - u_j)=M_1^{(1)}+1$ and $n(u_{i'}-u_{j'})\leq M_1^{(1)}$.
\end{itemize}
Under these consideration, let us define $A_0=\{U\in A|n(u_1 - u_2)=M_1^{(1)}+1\}$. It is obvious that $A=\cup_{i=0}^{p-1}\sigma^i(A_0)$ and $\cap_{i=1}^{p-1}\sigma^i(A_0)=\emptyset$. Therefore, $\gamma_{pp}(A)\leq p-1$.

On the other hand, we also need to bound $\gamma_{pp}\Big(\mathcal{A}_{\mathop{P}\limits^{\rightarrow},\bf{R}+{\bf e}}^=\Big)$ from above. Here, the set $\mathcal{A}_{\mathop{P}\limits^{\rightarrow},\bf{R}+{\bf e}}^=$ is defined as in (\ref{WorkingSet4}).
To this end, we point out that
\begin{itemize}
  \item we are working in $(H_{0,r}^2(\Omega))^N$ settings and therefore in $(C(\Omega))^N$;
  \item for any $U=(u_1,\cdots,u_N)\in \mathcal{A}_{\mathop{P}\limits^{\rightarrow},\bf{M}+{\bf e}}^=$, any couple $(u_i,u_j)$ of Case $(b,q)$ we have $n(u_i - u_j)=M_q^{(b)}$. Here, $b=1,\cdots,B$ and $q=1,\cdots,\frac{p-1}{2}$;
  \item for any $U=(u_1,\cdots,u_N)\in \mathcal{A}_{\mathop{P}\limits^{\rightarrow},\bf{M}+{\bf e}}^=$ and any $i,j=1,\cdots,N$ with $i\neq j$, $u_i\neq u_j$.
\end{itemize}
Write
$$A_1:=\mathcal{A}_{\mathop{P}\limits^{\rightarrow},\bf{M}+{\bf e}}^=\cap\Big\{U\in(H_{0,r}^2(\Omega))^N\Big|(u_1 - u_2)_1\geq0\Big\}$$
and
$$A_2:=\mathcal{A}_{\mathop{P}\limits^{\rightarrow},\bf{M}+{\bf e}}^=\cap\Big\{U\in(H_{0,r}^2(\Omega))^N\Big|(u_1 - u_2)_1\leq0\Big\}.$$
Here the function $(u_1 - u_2)_1$ denotes the first bump of $u_1 - u_2$. One can easily verify that
\begin{itemize}
  \item $\mathcal{A}_{\mathop{P}\limits^{\rightarrow},\bf{M}+{\bf e}}^= = \bigcup_{i=0}^{p-1}\sigma^i(A_1)$;
  \item $\cap_{i=0}^{p-1}\sigma^i(A_1)=\emptyset$.
\end{itemize}
By Assertion $(7)$ of Proposition \ref{prop:genus}, this is sufficient to imply $\gamma_{pp}\Big(\mathcal{A}_{\mathop{P}\limits^{\rightarrow},\bf{R}+{\bf e}}^=\Big)\leq p-1$. Therefore,
\begin{align}
\gamma_{pp}\Big(\mathcal{A}_{\mathop{P}\limits^{\rightarrow},\bf{M}+{\bf e}}^\leq\backslash\mathcal{D}\Big)\leq \gamma_{pp}\Big(\mathcal{A}_{\mathop{P}\limits^{\rightarrow},\bf{M}}^\leq\backslash\mathcal{D}\Big) +2(p-1).\nonumber
\end{align}
This is sufficient to imply the result.

\begin{flushright}
$\Box$
\end{flushright}

\subsection{More notations}\label{Subsection:Simplex}
In this subsection, we will continue to introduce the notations we will use.
We will briefly review certain settings in \cite{LiWang2021} with similar ones in \cite{IshiwataLi}.
\begin{itemize}
  \item Denote
  \begin{align}\label{def:H}
  H=\Bigg\{U= & (u_{1},\dots,u_{N})\in(H_{0,r}^{2}(\Omega))^{N}| n(u_{(b-1)p+j})\leq P_b\,\,\mbox{for}\, j=1,\dots,p,\,b=1,\cdots,B,\nonumber\\
   & n(u_{Bp+r})\leq Q_r\mbox{ for }r=1,\cdots,R\mbox{ and}\, \sum_{j=1}^N n(u_j) < p\sum_{j=1}^B P_j+\sum_{r=1}^R Q_r\Bigg\};
  \end{align}
\end{itemize}
We can now introduce the complete invariant set, which was used in \cite{LiuSun2001}.
\begin{itemize}
  \item We define
  \begin{align}\label{def:Cjq1}
  C_{(b-1)p+j,q}=\{U\in\mathcal{A}_{\mathop{P}\limits^{\rightarrow}}|\exists T\geq0\mbox{ s.t. }|\eta^T(U)_{(b-1)p+j,q}|_4\leq{\varepsilon}\}
  \end{align}
  for $q=1,\dots,P_{b}+1$ and $j=1,\dots,N$;
  \item We define
  \begin{align}\label{def:Cjq2}
  C_{Bp+r,q}=\{U\in\mathcal{A}_{\mathop{P}\limits^{\rightarrow}}|\exists T\geq0\mbox{ s.t. }|\eta^T(U)_{Bp+r,q}|_4\leq{\varepsilon}\}
  \end{align}
  for $r=1,\cdots,R$ and $q=1,\cdots,Q_r$.
\end{itemize}
Here, the set $\mathcal{A}_{\mathop{P}\limits^{\rightarrow}}$ is defined as in (\ref{WorkingSet1}) and the number $\varepsilon>0$ is a small number ensuring the validity of Proposition \ref{p:L4}.

Such notations are well-defined due to Subsection 2.3.
They present the constructions that naturally exist in the Sobolev space equipped with parabolic flow.
Let us now provide a brief interpretation of the meaning of these sets.
\begin{itemize}
  \item $H$ contains the vector-valued functions with less number of nodes than we prescribed;
  \item $C_{j,q}$ denotes the set of the functions whose $L^4$ norm of the $q$-th bump of the $j$-th component is less than $\varepsilon$.
\end{itemize}

Now we define some auxiliary functions.
We will use these functions to find initial data with certain properties. To begin with, let us recall the constructions in \cite{LiWang2021}. We start by dividing the domain $\Omega$.
\begin{itemize}
  \item [$(1).$] Divide the radial domain $\Omega$ into $B+R$ disjoint radial parts, say $\Omega_1^{(1)}$,..., $\Omega_B^{(1)}$; $\Omega_1^{(2)}$, ...,$\Omega_R^{(2)}$, ordered by the distance from the origin;
  \item [$(2).$] For the first group of the sub-domains, for any $b=1,\cdots,B$, divide $\Omega_b^{(1)}$ into $P_b+1$ radial parts. We denote them by $\Omega^{(1)}_{b,q}$ for $b=1,\cdots,B$ and $q=1,\cdots,P_b+1$;
  \item [$(3).$] For the second group of the sub-domains, for any $r=1,\cdots,R$, divide $\Omega_r^{(2)}$ into $Q_r+1$ radial parts. We denote them by $\Omega^{(2)}_{r,q}$ for $r=1,\cdots,R$ and $q=1,\cdots,Q_r+1$;
  \item [$(4).$] For each of the sub-domains $\Omega_{b,q}^{(1)}$'s and $\Omega_{r,q}^{(2)}$'s, we cut it into $K$ radial parts, written as $\Omega_{b,q,k}^{(1)}$'s and $\Omega_{r,q,k}^{(2)}$'s for $k=1,\cdots,K$.
\end{itemize}
Now we define the functions on them. For the first group of sets $\Omega_{b,q,k}^{(1)}$ with $b=1,\cdots,B$, $q=1,\cdots,P_b +1$ and $k=1,\cdots,K$,
\begin{itemize}
  \item [$(a).$] $w^{(1)}_{b,q,k}(t,x)=w^{(1)}_{b,q,k}(t,|x|)=w^{(1)}_{b,q,k}(t,r):\mathbb{S}^1\times\Omega^{(1)}_{b,q,k}\to[0,+\infty)$ of class $C^{4}$ and of compact support in $\mathbb{S}^1\times\Omega^{(1)}_{b,q,k}$;
  \item [$(b).$] $w^{(1)}_{b,q,k}(t,\cdot)\not\equiv0$ for any $t\in\mathbb{S}^1$ ;
  \item [$(c).$] $\mbox{supp}_x w^{(1)}_{b,q,k}(t,\cdot)\cap \mbox{supp}_x w^{(1)}_{b,q,k}\big(\frac{2\pi }{p}+t,\cdot\big)=\emptyset$ for any $t\in\mathbb{S}^1$.
\end{itemize}

\begin{remark}\label{r:support}
Different from \cite{LiWang2021}, we need to go into the details of such kind of auxiliary functions.
To properly compare their differences, the support of those functions must be located on $\mathbb{S}^1\times\Omega_{b,p,k}^{(1)}$.

We note that it is possible to let $n\big(w^{(1)}(t,\cdot)-w^{(1)}(t+\frac{2\pi q}{p})\big)\leq2$ for $q=1,\cdots,p-1$.
For illustrative purposes, we provide graphs displaying the supports of the auxiliary functions for the $p=2$ and $p=3$ cases on $\mathbb{S}^1\times\Omega_{b,p,k}^{(1)}$.
\begin{figure}[htbp]
\centering
\begin{minipage}[t]{0.48\textwidth}
\centering
\includegraphics[width=6cm]{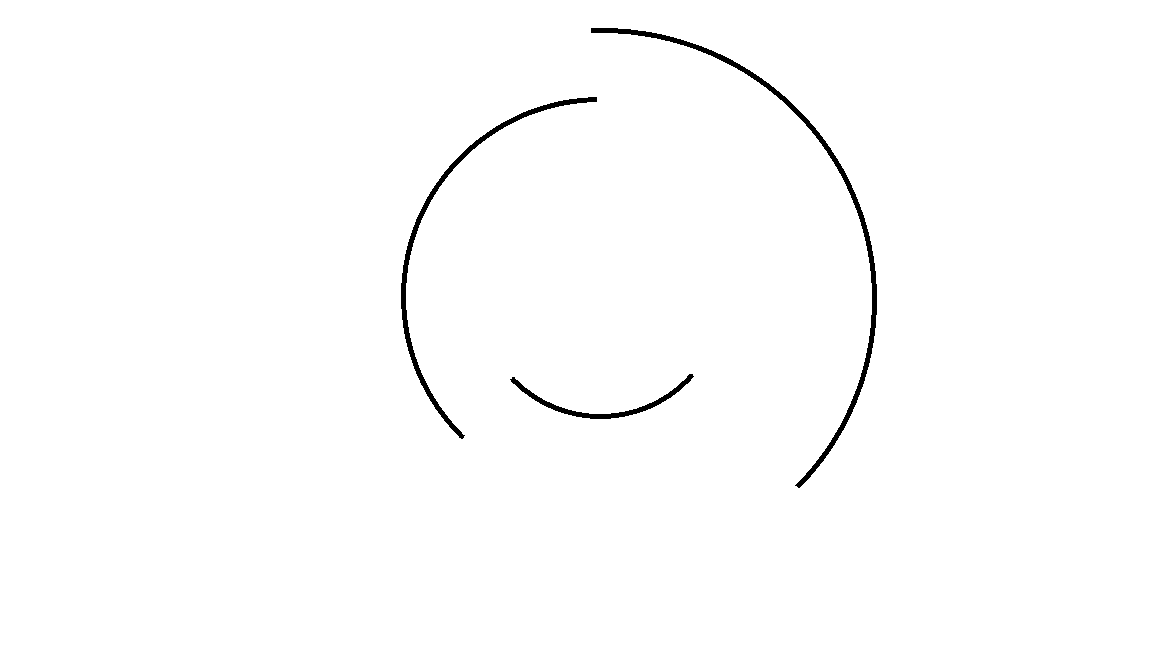}\label{FA}
\caption{The case for $p=2$}
\end{minipage}
\begin{minipage}[t]{0.48\textwidth}
\centering
\includegraphics[width=6cm]{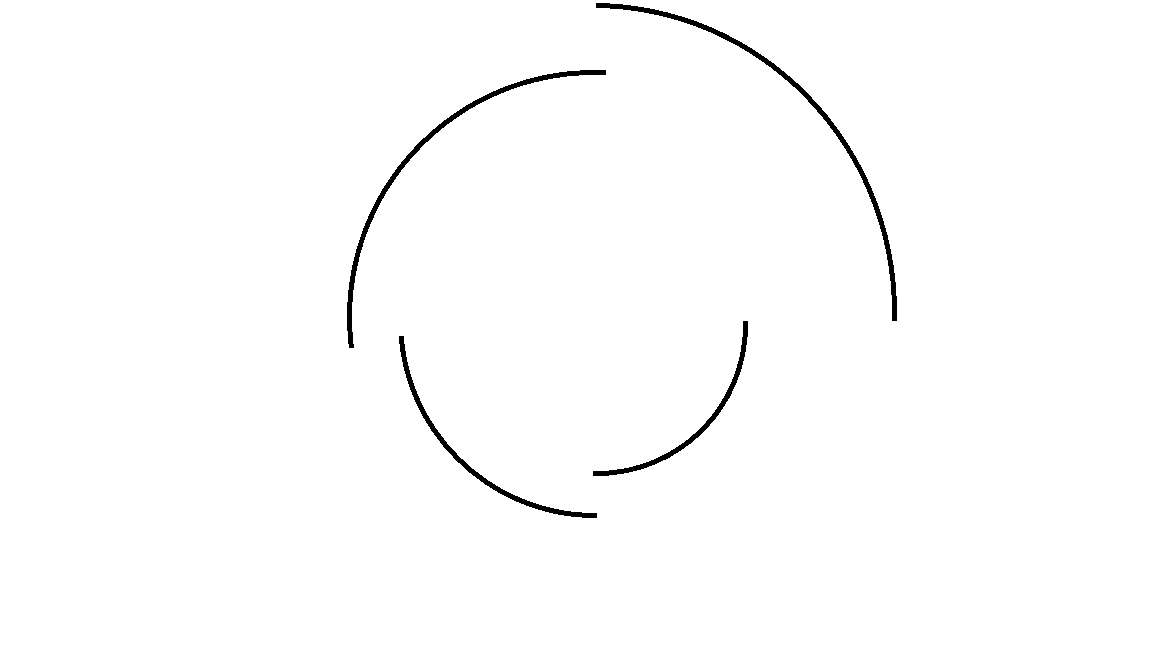}\label{FB}
\caption{The case for $p=3$}
\end{minipage}
\end{figure}

On $\mathbb{S}^1$-variable, we begin by dividing the domain equally into $p+1$ parts, with a slight overlaps at the endpoints since we need to ensure the non-triviality, i.e. Assertion $(b)$. This is shown in the figure. On each sub-arc, we select the location of supports and on each part build the auxiliary functions along them.

This construction guarantees that
$n\big(w^{(1)}(t,\cdot)-w^{(1)}(t+\frac{2\pi q}{p},\cdot)\big)\leq2$ for $q=1,\cdots,p-1$.
\end{remark}

Comparing to the auxiliary functions defined on $\mathbb{S}^1\times\Omega_{b,q,k}^{(1)}$'s, those defined on $\Omega_{r,q,k}^{(2)}$'s are simpler. It is sufficient to consider radial smooth functions $w_{r,q,k}^{(2)}:\Omega_{r,q,k,re}^{(2)}\to[0,+\infty)$.

Let us consider the complex space
$\mathbb{C}^{K(\sum_{b=1}^{B}P_b +\sum_{r=1}^R Q_r +B+R)}$. For $z:=(z_{b,q,k}^{(1)},z_{r,q,k}^{(2)})\in\mathbb{C}^{K(\sum_{b=1}^{B}P_b +\sum_{r=1}^R Q_r +B+R)}$,
\begin{align}
U_b(t,z)(x):=\sum_{q=1}^{P_b +1}\sum_{k=1}^K (-1)^{q+1}|z_{b,q,k}^{(1)}|\cdot w_{b,q,k}^{(1)}\Big(t+\mbox{arc}(z_{b,q,k}^{(1)}),x\Big),\nonumber
\end{align}
\begin{align}
V_r(z)(x):=\sum_{q=1}^{Q_r +1}\sum_{k=1}^K(-1)^{q+1}|z_{r,q,k}^{(2)}|\cdot w_{b,q,k}^{(2)}(x)\nonumber
\end{align}
and
\begin{align}
\psi(z):= \Bigg(U_1(0,z),\cdots,U_1\left(\frac{2\pi(p-1)}{p},z\right); & \cdots;U_B(0,z),\cdots,U_B\left(\frac{2\pi(p-1)}{p},z\right);\nonumber\\
& V_1(z),\cdots,V_R(x)\Bigg).\nonumber
\end{align}
Such a mapping satisfies $\psi:\mathbb{C}^{K(\sum_{b=1}^{B}P_b +\sum_{r=1}^R Q_r +B+R)}\to (H_{0,r}^2(\Omega))^N$ and $\psi(\mathbb{C}^{K(\sum_{b=1}^{B}P_b +\sum_{r=1}^R Q_r +B+R)})\subset(C_0^\infty(\Omega))^N$. Now we expand the simplex into the following form for the sake of computations.
\begin{align}\label{def:GK}
G_K=&\Bigg\{\Bigg(\sum_{k=1}^{K}\sum_{q=1}^{P_{1}+1}(-1)^{q+1}\alpha^{(1,1)}_{1,q,k}w^{(1)}_{1,q,k} (\theta^{(1)}_{1,q,k},x),\dots,\sum_{k=1}^{K}\sum_{q=1}^{P_{1}+1} (-1)^{q+1}\alpha_{1,q,k}^{(1,p)}w^{(1)}_{1,q,k}\Big(\frac{2\pi (p-1)}{p}+\theta^{(1)}_{1,q,k},x\Big);\nonumber\\
&\qquad\qquad\dots\dots\qquad\qquad\dots\dots\qquad\qquad\dots\dots\qquad\qquad\dots\dots;\nonumber\\
&\sum_{k=1}^{K}\sum_{q=1}^{P_{B}+1}(-1)^{q+1}\alpha^{(B,1)}_{1,q,k}w_{B,q,k}^{(1)}(\theta^{(1)}_{1,q,k},x), \dots,\sum_{k=1}^{K}\sum_{q=1}^{P_{B}+1} (-1)^{q+1}\alpha_{B,q,k}^{(B,p)}w_{B,q,k}^{(1)}\Big(\frac{2\pi (p-1)}{p}+\theta^{(1)}_{B,q,k},x\Big);\nonumber\\
&\sum_{k=1}^{K}\sum_{q=1}^{Q_1 +1}(-1)^{q+1}\alpha_{1,q,k}^{(2)} w_{1,q,k,}^{(2)}(x),\cdots\cdots,\sum_{k=1}^{K}\sum_{q=1}^{Q_R +1}(-1)^{q+1}\alpha_{R,q,k}^{(2)} w_{R,q,k}^{(2)}(x)\Bigg)\Bigg|\nonumber\\
&\alpha_{b,q,k}^{(1,j)},\alpha_{r,q,k}^{(2)}\geq0,\,\,\theta^{(1)}_{b,q,k}\in[0,2\pi),\,\,\mbox{for}\,\,\mbox{any}\,\,b=1,\dots,B,\,\,j=1\dots,p,\,\,q=1,\dots,P_{b}+1,\nonumber\\
&r=1,\cdots,R,\,k=1,\dots,K\Bigg\}.
\end{align}
The main difference between $G_K$ and $\psi(\mathbb{C}^{K(\sum_{b=1}^{B}P_b +\sum_{r=1}^R Q_r +B+R)})$ lies is that, in $G_K$ the component of dimensions are independent in each other.
In the following we construct a subspace in Euclidean space homeomorphic as $G_k$.
We begin by the notation of elements. Denote
\begin{align}
z^{(1)}:=(z_{pb+j-p,q,k}^{(1)})_{b=1,\cdots,B;j=1,\cdots,p; q=1,\cdots,P_b+1; k=1,\cdots,K}\nonumber
\end{align}
and
\begin{align}
z^{(2)}:=(z^{(2)}_{Bp+r,q,k})_{r=1,\cdots,R;q=1,\cdots,Q_r+1;k=1,\cdots,K}.\nonumber
\end{align}
Define
\begin{align}\label{def:X1}
X_1&=\Big\{z=(z^{(1)},z^{(2)})\in\mathbb{C}^{K(p\sum_{b=1}^B P_b+\sum_{r=1}^R Q_r+N)}\Big|\mbox{for any }b=1,\cdots,B;\,q=1,\cdots,P_b+1;\nonumber\\
&\quad k=1,\cdots,K,\mbox{ we have }\mbox{arc}(z_{(b-1)p+1,q,k}^{(1)})=\mbox{arc}(z_{(b-1)p+2,q,k}^{(1)})+ \frac{2\pi}{p}=\cdots\nonumber\\
&\quad=\mbox{arc}(z^{(1)}_{bp,q,k})+\frac{2\pi(p-1)}{p}\Big\}.
\end{align}
It is obvious that $X_1$ is homeomorphic to $G_K$.

On $X_1$, we define the following mapping which provides the partial standard rotation symmetry on the complex space
\begin{align}
\sigma_p(z^{(1)},z^{(2)})=(e^{\frac{2\pi i}{p}}z^{(1)},z^{(2)}).\nonumber
\end{align}
That is, $\sigma_p$ acts as rotations on the first component $z^{(1)}$, and $\sigma_p^p=Id$. The $\mathbb{Z}_p$-index $\gamma_{\sigma_p}$ generated by $\sigma_p$ can be defined in a more standard way. To be precise, we define it as follows.
\begin{itemize}
  \item $\mbox{Fix}_{\sigma_p}=\{U\in X_1|\sigma_p(U)=U\}$;
  \item $\mathcal{E}_{\sigma_p}=\{A\subset X_1\backslash\mbox{Fix}_{\sigma_p}|\sigma(A)=A\mbox{ and }A\mbox{ is compact}\}$;
  \item for any $A\in\mathcal{E}_{\sigma_p}$, define $\mathbb{I}_2(A):=\{m\in\mathbb{N}_+|\exists\mbox{ a continous }h:A\to\mathbb{C}^m\backslash\{0\}\mbox{ s.t. }h(\sigma_p z)=e^{\frac{2\pi i}{p}}h(z)\mbox{ }\forall U\in A\}$;
  \item for any $A\in\mathcal{E}_{\sigma_p}$, define
  \begin{equation}
\gamma_{\sigma_p}(A)=\left\{
\begin{aligned}
&\min\{m|m\in\mathbb{I}_2(A)\}&\mbox{ if }\mathbb{I}_2(A)\neq\emptyset, \nonumber\\
&\infty & \mbox{ if }\mathbb{I}_2(A)=\emptyset.\nonumber
\end{aligned}
\right.
\end{equation}
\end{itemize}
An analogue of Proposition \ref{prop:genus} for $\gamma_{\sigma_p}$ holds.
\begin{proposition}\label{prop:genus0} For any $A,B\in\mathcal{E}_{\sigma_p}$, we have
\begin{itemize}
        \item [$(1).$] If $A\subset B$, then $\gamma_{\sigma_p}(A)\leq\gamma_{\sigma_p}(B)$;
        \item [$(2).$] $\gamma_{\sigma_p}(A\cup B)\leq\gamma_{\sigma_p}(A)+\gamma_{\sigma_p}(B)$;
        \item [$(3).$] if $g:A\to X_1\backslash \mbox{Fix}_{\sigma_p}$ is continuous and satisfies
        $g(\sigma_p(z))=\sigma_p g(z)$ for all $ z\in A$,
        then $$\gamma_{\sigma_p}(A)\leq\gamma_{\sigma_p}(\overline{g(A)});$$
        \item [$(4).$] if $\gamma_{\sigma_p}(A)>1$, then $A$ is an infinite set;
        \item [$(5).$] if $A$ is compact and $\gamma_{\sigma_p}(A)<\infty$, then there exist an open $\sigma_p$-invariant neighbourhood $\mathcal{N}$ of $A$
        such that $\gamma_{\sigma_p}(A)=\gamma_{\sigma_p}(\overline{\mathcal{N}})$;
        \item [$(6).$] if $S$ is the boundary of a bounded neighbourhood of the origin in a $m$-dimensional complex linear space such that
       $e^{\frac{2\pi i}{p}}U \in S $ for any $U\in S$, and $\Psi:S\to X_1\backslash \mbox{Fix}_{\sigma_p} $ is continuous and satisfies for any $U\in S$,
       $\Psi(e^{\frac{2\pi i}{p}}U)=\sigma(\Psi(U))$, then $\gamma_{\sigma_p}(\Psi(S))\geq m$;
       \item [$(7).$] Let $A$ be a closed set such that $A\subset X_1\setminus \mbox{Fix}_{\sigma_p}$, and $\cap_{i=0}^{p-1} \sigma^i(A) =\emptyset $.
Then $\gamma_{\sigma_p}(\mathbb \mathbb{Z}_p(A))\leq p-1$.
   \end{itemize}
\end{proposition}

\subsection{The sets in which the nodes of the components are not vanishing}\label{Subsection:nonvanishing}

In this part, we combine the calculations on genus in \cite{LiuLiuWang2015} and \cite{LiWang2021}.
Our idea can be summarized as follows.
We like to find a linking $G_K\cap\partial\mathcal{A}\backslash\Big(\bigcup_{j,q}C_{j,q}\cup H\Big)$ on $\partial\mathcal{A}$ that admits the invariance. However, the genus of this linking is hard to compute.
Instead, we will find an auxiliary set $h(G_K\cap\overline{\mathcal{A}})\cap F_{2\varepsilon}$, whose genus can be computed. The genus $\gamma_{pp}(G_K\cap\partial\mathcal{A}\backslash(\bigcup_{j,q}C_{j,q}\cup H))$ is obtained by a homeomorphic argument.
By a homeomorphism induced by parabolic flow, we obtain that the genus of $h(G_K\cap\overline{\mathcal{A}})\cap F_{2\varepsilon}$ is a lower bound of the genus of $G_K\cap\partial\mathcal{A}\backslash\Big(\bigcup_{j,q}C_{j,q}\cup H\Big)$ on $\partial\mathcal{A}$.
Therefore, in order to find a lower estimate on
$$\gamma_{pp}\Bigg(G_K\cap\partial\mathcal{A}\backslash\Big(\bigcup_{j,q}C_{j,q}\cup H\Big)\Bigg),$$
we need to
\vs
\begin{itemize}
  \item [$\mbox{Step 1}.$] Show the relation between $G_K\cap\partial\mathcal{A}\backslash\Big(\bigcup_{j,q}C_{j,q}\cup H\Big)$ and $h(G_K\cap\overline{\mathcal{A}})\cap F_{2\varepsilon}$ in $\mathcal{A}$;
  \item [$\mbox{Step 2}.$] Compute a lower bound of $\gamma_{pp}\big(h(G_K\cap\overline{\mathcal{A}})\cap F_{2\varepsilon}\big)$ in $\mathcal{A}$.
\end{itemize}

\vs

Denote
\begin{align}\label{definitionDk}
D_K := G_K \cap \partial \mathcal{A} \backslash \Big( \bigcup_{j,q} C_{j,q} \cup H \Big)
\end{align}
and
$$T^{\varepsilon}(U):=\inf\Big\{T\geq0\Big|\exists(j,q)\mbox{ admissible such that }|\eta^T(U)_{j,q}|_4\leq\rho\mbox{ or }\eta^T(U)\in H \Big\}.$$
\begin{claim}
$T^\varepsilon(U)$ is continuous in $U\in\mathcal{A}_{\mathop{P}\limits^{\rightarrow}}\cap G$. Here, the sets $\mathcal{A}_{\mathop{P}\limits^{\rightarrow}}$ and $G_K$ are defined as in (\ref{WorkingSet1}) and in (\ref{def:GK}), respectively.
\end{claim}
We refer \cite[Lemma 3.4]{LiWang2021} for the proof.

Now we select a small positive number $\varepsilon>0$ and consider the following cut-off function.
$$\phi^{\varepsilon}(U):=\frac{d_{(H^1)^N}(U,G_K\cap\overline{\mathcal{ A}}\backslash(D_K)_{2\varepsilon})}{d_{(H^1)^N}(U,G_K\cap\overline{\mathcal{A}}\backslash(D_K)_{2\varepsilon}) +d_{(H^1)^N}(U,(D_K)_{\varepsilon})}.$$
Here, $d_{(H^1)^N}(\cdot,\cdot)$ denotes the distance function in $(H^1_{0,r}(\Omega))^N$. It is evident that $\phi^{\varepsilon}(U)$ is locally Lipschitzian on $G\cap\overline{\mathcal{A}}$.
Define
\begin{align}\label{def:h}
h(U):=\eta^{T^\varepsilon(U) \cdot \phi^{\varepsilon}(U)}(U)
\end{align}
and
\begin{align}\label{def:F}
F_{2\varepsilon}:=\left\{U\in\mathcal{A}_{\mathop{P}\limits^{\rightarrow}}\big||u_{j,q}|_4=2\varepsilon\right\}.
\end{align}
We will compute a lower bound for $\gamma_{pp}(h(G_K\cap\overline{\mathcal{A}})\cap F_{2\varepsilon})$.
To this end, we first study the mapping $h$.

\begin{lemma}
The mapping $h:\overline{\mathcal{A}}\cap G_K\to h(\overline{\mathcal{A}}\cap G_K)$ is a homeomorphism. Here, the sets $\mathcal{A}$ and $G_K$ are defined as in (\ref{def:A}) and in (\ref{def:GK}), respectively.
\end{lemma}
\noindent{\bf Proof.}
We only need to check that the mapping $h$ is a bijection. The rest of the proof follows immediately from \cite[Theorem 7.8/pp. 19]{Bredon1993}, the continuity of $h$ and the compactness of the set $\overline{\mathcal{A}}\cap G_K$. We argue it by contradiction. Let us assume that for two different $U_1$ and $U_2$ on $\overline{\mathcal{A}}\cap G_K$, $h(U_1)=h(U_2)$. Due to the definition of the mapping $h$, we divide the deduction into the following two cases: (1). the functions $U_1$ and $U_2$ are on the same flow line; (2). otherwise.

{\bf (1). There is a $t_0\in(0,\phi_1(U_1)\cdot T^{\varepsilon'}(U_1)]$ such that $U_2=\eta^{t_0}(U_1)$.}
Denote the "inverse" flow line $\theta^t(U_2)=\theta^t\left(\eta^{t_0}(U_1)\right):=\eta^{t_0-t}(U)$ for $t\in[0,t_0]$. It is easy to see that the inverse flow line satisfies the following initial value problem
\begin{equation}
    \left\{
   \begin{array}{lr}
     -\frac{\partial}{\partial t}v_{j}-{\Delta}v_{j}+ v_{j}=v^{3}_{j}+\beta\sum_{i\neq j} v_{j}v_{i}^{2} \mbox{ in }\Omega ,\nonumber\\
     (v_1(0,x),\cdots,v_N(0,x))=\eta^{t_2}(U).\nonumber
   \end{array}
   \right.
\end{equation}
Notice that both $\theta^0(U_2)=U_2=\eta^{t_0}(U_1)$ and $\theta^{t_0}(U_2)=\eta^{0}(U_1)=U_1$ are compactly supported. According to Theorem \ref{t:HeatExist} and Theorem \ref{t:backward}, $\eta^{t}(U_1)\equiv0$ for $t\in(0,t_0)$. This contradicts the continuity of the flow. The above method is also valid for the case when both of the vectors $U_1$ and $U_2$ have no trivial component.

{\bf (2). For $i=1,2$, for any $t\in(0,\phi_1(U_i)\cdot T^{\varepsilon}(U_i)]$, $U_{i^*}\neq\eta^{t}(U_i)$.} Here, $i^*=2$ when $i=1$ and $i^*=1$ when $i=2$. We introduction the notations:
$$t_1=\inf\{t>0|\eta^t(U_1)\cap\left(\cup_{s>0}\eta^s(U_2)\right)\neq\emptyset\};$$
$$t_2=\inf\{t>0|\eta^t(U_2)\cap\left(\cup_{s>0}\eta^s(U_1)\right)\neq\emptyset\}.$$
Note that in this case, both of the vectors $U_1$ and $U_2$ have no trivial components. Otherwise, we have $T^{\varepsilon}(U_1)=0$. Furthermore, we have $t_1,t_2\in(0,+\infty)$. Now we divide the discussion into two cases:

{\bf (2.1). The case of $t_1=t_2$.} For the sake of convenience, we denote this number by $t_1$. In this case, $\eta^{t_1}(U_1)=\eta^{t_1}(U_2)$. Denote the function $U_3(t)=\eta^{t_1-t}(U_1)-\eta^{t_1-t}(U_2)$ for $t\in[0,t_1]$. Due to Corollary \ref{coro:LinftyBdd}, notice that there is a constant $C>0$, for any $j=1,\cdots,N$ and any $(x,t)\in\Omega\times[0,t_1]$,
$$|\partial_t u_{3,j}+\Delta u_{3,j}|\leq C|u_{3,j}|.$$
Here, $u_{3,j}$ is the $j$-th component of $U_3$. Thus, we have a contradiction with the help of Theorem \ref{t:backward}, $U_3(0)=0$ and $U_3(t_1)\neq0$.

{\bf (2.2). The case of $t_1\neq t_2$.} Without loss of generality, let us assume that $t_1>t_2$. Denote $U'_1=\eta^{t_1-t_2}(U_2)$. Due to Theorem \ref{t:backward}, we note that the support of each component of $U'_1$ is $\Omega$. Applying a similar procedure as in {\bf Case (2.1)}, we will have a contradiction with $U'_1-U_2\neq0$. The latter is obvious since $U_2$ is compactly supported.

In summary, the mapping $h$ is 1-1. The claim is proved.

\begin{flushright}
$\Box$
\end{flushright}


The next lemma follows immediately.

\begin{lemma}\label{l:GenusLowerBdd}
It holds for small $\varepsilon>0$ that
$$\gamma_{pp}(D_K)=\gamma_{pp}\Bigg(G_K\cap\partial\mathcal{A}\backslash\Big(\bigcup_{j,q}C_{j,q}\cup H\Big)\Bigg)\geq\gamma_{pp}(h(G_K\cap\overline{\mathcal{A}})\cap F_{2\varepsilon}).$$
Here, the index $\gamma_{pp}$, the sets $D_K$, $G_K$, $\mathcal{A}$, $C_{j,q}$, $H$ and $F_{2\varepsilon}$ and the map $h$ are defined as in (\ref{def:GammaPP}), (\ref{definitionDk}), (\ref{def:GK}), (\ref{def:A}), (\ref{def:Cjq1}-\ref{def:Cjq2}), (\ref{def:H}), (\ref{def:F}) and (\ref{def:h}), respectively.
\end{lemma}

\noindent{\bf Proof.}
The result holds evidently when $h(G_K\cap\overline{\mathcal{A}})\cap F_{2\varepsilon}=\emptyset$. When $h(G_K\cap\overline{\mathcal{A}})\cap F_{2\varepsilon}\neq\emptyset$, notice that $h$ is a homeomorphism, which implies that the notation $h^{-1}\left(h(G_K\cap\overline{\mathcal{A}})\cap F_{2\varepsilon}\right)$ is well-defined. Now we locate $h^{-1}\left(h(G_K\cap\overline{\mathcal{A}})\cap F_{2\varepsilon}\right)$ in $G_K$.
\begin{claim}\label{c:emptyset}
For sufficiently small $\varepsilon'>0$, we get
$$h^{-1}\left(h(G_K\cap\overline{\mathcal{A}})\cap F_{2\varepsilon}\right)\cap\Big(G\cap\overline{\mathcal{A}}\backslash (D_K)_{\varepsilon'}\Big)=\emptyset.$$
Here, $(D_K)_{\varepsilon'}$ is the $\varepsilon'$-neighbourhood of the set $D_K$.
\end{claim}
If this holds, we get
$$h^{-1}\left(h(G_K\cap\overline{\mathcal{A}})\cap F_{2\varepsilon}\right)\subset G_K\cap\overline{\mathcal{A}}\cap (D_K)_{\varepsilon'}.$$
This gives that
$$h^{-1}\left(h(G_K\cap\overline{\mathcal{A}})\cap F_{2\varepsilon}\right)\subset \overline{(D_K)_{\varepsilon'}}.$$
For small $\varepsilon'>0$, it holds that
\begin{align}
\gamma_{pp}\left(h(G_K\cap\overline{\mathcal{A}})\cap F_{2\varepsilon}\right)=&\gamma_{pp}\left(h^{-1}\left(h(G_K\cap\overline{\mathcal{A}})\cap F_{2\varepsilon}\right)\right)\leq\gamma_{pp}\left(\overline{(D_K)_{\varepsilon'}}\right) =\gamma_{pp}(D_K)\nonumber\\
=&\gamma_{pp}\Bigg(G_K\cap\partial\mathcal{A}\backslash\Big(\bigcup_{j,q}C_{j,q}\cup H\Big)\Bigg)\nonumber
\end{align}
This proves Lemma \ref{l:GenusLowerBdd}.

Now we prove Claim \ref{c:emptyset}.
Assuming that $h^{-1}\left(h(G_K\cap\overline{\mathcal{A}})\cap F_{2\varepsilon}\right)\cap\Big(G_k\cap\overline{\mathcal{A}}\backslash (D_K)_{\varepsilon'}\Big)\neq\emptyset$,
we select a point $U_0\in h^{-1}\left(h(G_K\cap\overline{\mathcal{A}})\cap F_{2\varepsilon}\right)\cap\Big(G_k\cap\overline{\mathcal{A}}\backslash (D_K)_{\varepsilon'}\Big)$.
On one hand, $U_0\in h^{-1}\left(h(G_K\cap\overline{\mathcal{A}})\cap F_{2\varepsilon}\right)$. Then, $h(U_0)\in h(G_K\cap\overline{\mathcal{A}})\cap F_{2\varepsilon}$. For any $(j,q)$ admissible, we have that $|h(U_0)_{j,q}|_4=2\varepsilon$.
On the other hand, $U_0\in G_k\cap\overline{\mathcal{A}}\backslash (D_K)_{\varepsilon'}$.
Let $\varepsilon'$ be a positive number such that
\begin{align}
\phi^\varepsilon|_{G_K\cap\overline{\mathcal{A}}\backslash(D_K)_{\varepsilon'}}=Id.\nonumber
\end{align}
Hence,
\begin{align}
h|_{G\cap\overline{\mathcal{A}}\backslash (D_K)_{\varepsilon'}}(\cdot)=\eta^{T^{\varepsilon}(\cdot)}(\cdot).\nonumber
\end{align}
Then, we get $h(U_0)=\eta^{T^{\varepsilon}(U_0)}(U_0)$. Then there exists a $(j_0,q_0)$ admissible such that $|h(U_0)_{j_0,q_0}|_4\leq\varepsilon$ or $h(U_0)\in H$. This is a contradiction. Therefore, Claim \ref{c:emptyset} holds and Lemma \ref{l:GenusLowerBdd} follows.

\begin{flushright}
$\Box$
\end{flushright}

\begin{lemma}
For any large $K>0$ and small $\varepsilon>0$, $$\gamma_{pp}(h(G_K\cap\overline{\mathcal{A}})\cap F_{2\varepsilon})\geq K-(p-1)\Big(\sum_{b=1}^B P_b +B\Big).$$
Here, the index $\gamma_{pp}$, the map $h$, the sets $G_K$, $\mathcal{A}$ and $F_{2\varepsilon}$ are defined as in (\ref{def:GammaPP}), (\ref{def:h}), (\ref{def:GK}), (\ref{def:A}) and (\ref{def:F}), respectively. The numbers $P_b$, $B$ and $p$ are prescribed in Theorem \ref{t:existence}.
\end{lemma}
\noindent{\bf Proof.}
To begin with, we introduce a few notations.
\begin{itemize}
  \item $\mathcal{O}=\{z\in X_1| h(z)\in\overline{M_{2\varepsilon}}\}$;
  \item $M_{2\varepsilon}=\{U\in \mathcal{A}_{\mathop{P}\limits^{\rightarrow}}||u_{j,q}|_4\leq2\varepsilon\mbox{ for any admissible }(j,q)\}$.
\end{itemize}
Recall that we define
\begin{align}
\mathcal{A}_{\mathop{P}\limits^{\rightarrow}}:=\Big\{U\in\overline{\mathcal{A}}\cap(H_{0,r}^2(\Omega))^N \Big|n(u_{(b-1)p+j})=P_b, & \,\,b=1,\cdots,B; \,j=1\cdots,p\nonumber\\
&\mbox{ and }n(u_{pB+r})=Q_i,\,\,r=1,\cdots,R\Big\}.\nonumber
\end{align}
Thanks to Borsuk's theorem,
\begin{align}\label{ineq:genus01}
\gamma_{\sigma_p}(\partial\mathcal{O})\geq K\Big(\sum_{b=1}^{B}P_b +\sum_{r=1}^R Q_r +B+R\Big).
\end{align}
since $\psi(\mathbb{C}^{K(\sum_{b=1}^{B}P_b +\sum_{r=1}^R Q_r +B+R)})\subset G_k$. Here, $\gamma_{\sigma_p}$ is the $\mathbb{Z}_p$-genus generated by the action
\begin{align}
\sigma_p:\mathbb{C}^{K(\sum_{b=1}^{B}P_b +\sum_{r=1}^R Q_r +B+R)}&\to\mathbb{C}^{K(\sum_{b=1}^{B}P_b +\sum_{r=1}^R Q_r +B+R)}\nonumber\\
(z_1,\cdots,z_{K(\sum_{b=1}^{B}P_b +\sum_{r=1}^R Q_r +B+R)})&\mapsto(e^{\frac{2\pi i}{p}}z_1,\cdots,e^{\frac{2\pi i}{p}}z_{K(\sum_{b=1}^{B}P_b +\sum_{r=1}^R Q_r +B+R)}).\nonumber
\end{align}
Here, the numbers $B$, $P_b$, $R$, $Q_r$ are prescribed in Theorem \ref{t:existence}.
By a routine computation as in \cite[Lemma 4.8]{LiWang2021}, we can conclude that
\begin{align}\label{ineq:genus001}
\gamma_{\sigma_p}(\partial\mathcal{O})\leq K\Big(\sum_{b=1}^{B}P_b +\sum_{r=1}^R Q_r +B+R\Big).
\end{align}
We prove it in Appendix.
(\ref{ineq:genus01}) and (\ref{ineq:genus001}) together give
\begin{align}\label{ineq:genus1}
\gamma_{\sigma_p}(\partial\mathcal{O})= K\Big(\sum_{b=1}^{B}P_b +\sum_{r=1}^R Q_r +B+R\Big).
\end{align}

Denote
\begin{itemize}
  \item $R_{(b-1)p+j,q}^{2\varepsilon}=\{z\in X_1||h(z)_{(b-1)p+j,q}|_4\leq{2\varepsilon}\}$ for $q=1,\dots,P_{b}+1$ and $j=1,\dots,N$;
  \item $R_{Bp+r,q}^{2\varepsilon}=\{z\in X_1||h(z)_{pB+r,q}|_4\leq{2\varepsilon}\}$ for $r=1,\cdots,R$ and $q=1,\cdots,Q_r$;
  \item $S_{(b-1)p+j,q}^{2\varepsilon}=\{z\in X_1||h(z)_{(b-1)p+j,q}|_4 ={2\varepsilon}\}$ for $q=1,\dots,P_{b}+1$ and $j=1,\dots,N$;
  \item $S_{Bp+r,q}^{2\varepsilon}=\{z\in X_1||h(z)_{Bp+r,q}|_4=2\varepsilon\}$ for $r=1,\cdots,R$ and $q=1,\cdots,Q_r$.
\end{itemize}
Let
\begin{align}\label{ineq:genus2}
L:=\partial\mathcal{O}\backslash(A_1\cup A_2 \cup A_3).
\end{align}
Here,
\begin{itemize}
  \item $A_1=\{z\in X_1|\exists b_U=1,\cdots,B,\exists q_U=1,\cdots, P_{b_U}+1,\exists i_U,i'_{U}=1,\cdots,p\mbox{ with }i_U\neq i'_U\mbox{ such that }z\in R_{(b_U -1)p + i_U , q_U}^{2\varepsilon}\cap S_{(b_U -1)p +i'_U , q_U}^{2\varepsilon}\}$;
  \item $A_2=\cup_{b=1}^B\cup_{q=1}^{P_b +1}\cap_{j=1}^p R_{(b-1)p + j,q}^{2\varepsilon}$;
  \item $A_3=\cup_{r=1}^R\cup_{q=1}^{Q_r + 1} R_{Bp+r,q}^{2\varepsilon}$.
\end{itemize}
Then, $h(L)=h(G_K\cap\overline{\mathcal{A}})\cap F_{2\varepsilon}$.
Then, it is sufficient to consider the genus of $\partial\mathcal{O}\backslash(A_1\cup A_2 \cup A_3)$.
On the other hand, by the computation methods in \cite{LiWang2021}, it is easy to verify that
\begin{align}\label{ineq:genus11}
  \gamma_{\sigma_p}(A_1)\leq(p-1)\Big(\sum_{b=1}^B P_b +B\Big);
\end{align}
and
\begin{align}\label{ineq:genus22}
  \gamma_{\sigma_p}\big((A_2 \cup A_3)\cap\partial\mathcal{O}\backslash A_1\big)\leq K\Big(\sum_{b=1}^{B}P_b +\sum_{r=1}^R Q_r +B+R-1\Big).
\end{align}
Since the proofs of (\ref{ineq:genus11}) and (\ref{ineq:genus22}) are routine but long, we leave them to the appendix.
Therefore, (\ref{ineq:genus1}) and (\ref{ineq:genus2}) give $\gamma_{\sigma_p}(L)\geq K-(p-1)\Big(\sum_{b=1}^B P_b +B\Big)$.
Applying Lemma \ref{l:2genus}, $\gamma_{pp}(h(G_K\cap\overline{\mathcal{A}})\cap F_{2\varepsilon})\geq K-(p-1)\Big(\sum_{b=1}^B P_b +B\Big)$.
\begin{flushright}
$\Box$
\end{flushright}

We summarize the computation in this subsection into the following claim.
\begin{corollary}\label{coro:lowerDk}
It holds that $$\gamma_{pp}\Bigg(G_K\cap\partial\mathcal{A}\backslash\Big(\bigcup_{j,q}C_{j,q}\cup H\Big)\Bigg)\geq K-(p-1)\Big(\sum_{b=1}^B P_b +B\Big).$$
Here, the index $\gamma_{pp}$, the sets $G_K$, $\mathcal{A}$ and $H$ are defined as in (\ref{def:GammaPP}), (\ref{def:GK}), (\ref{def:A}) and (\ref{def:H}), respectively. The numbers $P_b$, $B$ and $p$ are prescribed in Theorem \ref{t:existence}.
\end{corollary}

\begin{remark}\label{r:DDempty}
In closing this subsection, we note that
\begin{align}
D_K\cap\mathcal{D}=\emptyset.
\end{align}
Here, the set $D_K$ is defined as in (\ref{definitionDk}) and $\mathcal{D}$ in (\ref{Set:D}).
\end{remark}

\subsection{The sets starting from which the comparisons hold on the flow lines}\label{Subsection:Comparison}

\vskip .2in

In the last subsection, we obtain that
\begin{align}\label{ineq:model}
K-(p-1)\sum_{b=1}^B (P_b +1)\leq\gamma_{pp}(D_K)\leq K\Big(\sum_{b=1}^{B}(P_b+1) +\sum_{r=1}^R (Q_r+1)\Big).
\end{align}
Recall that we assumed in (\ref{def:Ks}) that
\begin{align}\label{ASSUMPTION2}
K_{s+1}-8\sum_{b=1}^B (P_b+1)\cdot K_s(p-1)^2-5B(p-1)^2-\sum_{b=1}^B (P_b +1)(p-1)=1
\end{align}
with $K_1=8\sum_{b=1}^B (P_b +1)p(p+1)+5Bp^2$. Select the sequence of sets $\{D_{K_s}\}_{s=1}^\infty$.
It is known that
\begin{itemize}
  \item [$(1).$] $K_{s+1} - (p-1)\sum_{b=1}^B (P_b +1)\leq\gamma_{pp}(D_{K_{s+1}})$;
  \item [$(2).$] $D_{K_{s+1}}\cap\mathcal{D}=\emptyset$;
  \item [$(3).$] $\gamma_{pp}\Big(\mathcal{A}_{\mathop{P}\limits^{\rightarrow},\bf{M}_s}^{\leq} \backslash\mathcal{D}\Big)\leq (p-1)\left(2\inf_{b=1,...,B;q=1,...,\frac{p-1}{2}}M_{q,s}^{(b)}+3\right)$
\end{itemize}
with the sequence of matrices
\begin{align}
{\bf M}_s=\left(
  \begin{array}{cccc}
    M_{1,s}^{(1)} & M_{2,s}^{(1)} & \cdots & M_{\frac{p-1}{2},s}^{(1)} \\
    M_{1,s}^{(2)} & M_{2,s}^{(2)} & \cdots & M_{\frac{p-1}{2},s}^{(2)} \\
    \vdots & \vdots & \ddots & \vdots \\
    M_{1,s}^{(B)} & M_{2,s}^{(B)} & \cdots & M_{\frac{p-1}{2},s}^{(B)} \\
  \end{array}
\right).
\end{align}
also to be settled. Here, Assertion $(1)$ is due to Corollary \ref{coro:lowerDk}. Assertions $(2)$ and $(3)$ are ensured by Remark \ref{r:DDempty} and Lemma \ref{l:CompareSet}, respectively. The $(b,q)$-element of the matrix $R_{q,s}^{(b)}$ describes the maximum of the nodal numbers of the comparison of the couples in Case $(b,q)$, cf. Subsection \ref{Subsection:workingset}.
To continue, we first analyze $\eta^t(D_{K_{s+1}})=\{\eta^t(U)|U\in D_{K_{s+1}}\}$ for $t>0$. Notice that
\begin{itemize}
  \item $D_{K_{s+1}}\subset G_{K_{s+1}}$;
  \item for any $b=1,\cdots,B$ and any $i,i'=1,\cdots,p$ and $i\neq i'$, Remark \ref{r:support} implies that
      \begin{align}\label{ineq:ComparisonUpper}
      n\big(U(0)_{(b-1)p+i}-U(0)_{(b-1)p+i'}\big)\leq 4 (P_b+1)\cdot K_{s+1} +1.
      \end{align}
\end{itemize}
Here, the sets $D_K$ and $G_K$ are defined as in (\ref{definitionDk}) and (\ref{def:GK}), respectively.
Applying Corollary \ref{c:nodalnumber}, for any $t>0$, $n\big(U(t)_{(b-1)p+i}-U(t)_{(b-1)p+i'}\big)\leq 4 P_b\cdot K_{s+1} +1$ under the above notations. Letting $R_{q,s}^{(b)}=4P_b\cdot K_{s+1}+1$ for any $b=1,\cdots,B$, $q=1,\cdots,\frac{p-1}{2}$ and $s=1,2,\cdots$, we get
\begin{align}\label{ineq:Upper}
\gamma_{pp}(D_{K_{s+1}})\leq (p-1)\Big(8 K_{s+1}\min_b (P_b+1) +5\Big).
\end{align}
This is Lemma \ref{l:CompareSet}.
Now let us consider the sequence of matrices
\begin{align}
{\bf M}_{s}=\left(
  \begin{array}{cccc}
    4(P_1+1)\cdot K_{s+1} +1 & 4(P_1+1)\cdot K_{s+1} +1 & \cdots & 4(P_1+1)\cdot K_{s+1} +1 \\
    4(P_2+1)\cdot K_{s+1} +1 & 4(P_2+1)\cdot K_{s+1} +1 & \cdots & 4(P_2+1)\cdot K_{s+1} +1 \\
    \cdots & \cdots & \cdots & \cdots \\
    4(P_B+1)\cdot K_{s+1} +1 & 4(P_B+1)\cdot K_{s+1} +1 & \cdots & 4(P_B+1)\cdot K_{s+1} +1
  \end{array}
\right),
\end{align}
for $b=1,\cdots,B$ and $s=1,2,\cdots$. Define ${\bf M}_{s,b,j}$ to be the matrix ${\bf M}_s$ with its $(b,j)$-element replaced by $4(P_b+1)\cdot K_s +1$. By Lemma \ref{l:CompareSet}, we get
\begin{align}\label{ineq:GenusA}
\gamma_{pp}\Big(\mathcal{A}_{\mathop{P}\limits^{\rightarrow},\bf{M}_{s,b,j}}^{\leq}\backslash\mathcal{D} \Big)\leq(p-1)\big[8(P_b+1)K_s +5\big].
\end{align}
Define the set
\begin{align}\label{SetE}
E_{n,s}:=\overline{\Bigg\{U\in D_{K_{s+1}}\Bigg| \eta^n(U)\notin\bigcup_{b=1}^B \bigcup_{j=1}^{\frac{p-1}{2}}\mathcal{A}_{\mathop{P}\limits^{\rightarrow},\bf{M}_{s,b,j}}^{\leq}\Bigg\}}.
\end{align}
Such a set contains the elements $U$ in $D_{K_{s+1}}$ with $\eta^n(U)$ satisfies that for any $b=1,\cdots,B$ and any $i_1, i_2=1,\cdots,p$ with $i_1\neq i_2$, we have $n(\eta^n(U)_{(b-1)p+i_1}-\eta^n(U)_{(b-1)p+i_2})\leq 4(P_b+1)\cdot K_{s+1}$+1.
By (\ref{ineq:GenusA}), we get
\begin{align}\label{ineq:Upper1}
\gamma_{pp}\Bigg(\bigcup_{b=1}^B\bigcup_{j=1}^{\frac{p-1}{2}} \mathcal{A}_{\mathop{P}\limits^{\rightarrow},\bf{M}_{s,b,j}}^{\leq}\backslash\mathcal{D} \Bigg)\leq 8(p-1)^2\sum_{b=1}^B (P_b+1)\cdot K_s+5B(p-1)^2.
\end{align}
To ensure this, it is sufficient to assume that
\begin{align}\label{ASSUMPTION1}
\max_b (P_b+1)\cdot K_s\leq\min_b (P_b+1)\cdot K_{s+1}.
\end{align}
Then, the following claim is evident.
\begin{claim}
It holds that
$$\gamma_{pp}\big(E_{n,s}\big)\geq K_{s+1}-8(p-1)^2\sum_{b=1}^B (P_b+1)\cdot K_s-5B(p-1)^2-(p-1)\sum_{b=1}^B (P_b +1)$$
whenever the lower bound is positive.
\end{claim}
\noindent{\bf Proof.}
By a direct computation, we get
\begin{align}
K_{s+1} - (p-1)\sum_{b=1}^B (P_b +1)&\leq\gamma_{pp}(D_{K_{s+1}})\leq\gamma_{pp}(E_{n,s})+\gamma_{pp}(D_{K_{s+1}}\backslash E_{n,s})\nonumber\\
&\leq\gamma_{pp}(E_{n,s})+\gamma_{pp}\big(\eta^n(D_{K_{s+1}}\backslash E_{n,s})\big)\nonumber\\
&\leq\gamma_{pp}(E_{n,s})+8(p-1)^2\sum_{b=1}^B (P_b+1)\cdot K_s+5B(p-1)^2.\nonumber
\end{align}
Here, in the last inequality, we use (\ref{ineq:Upper1}).
Then,
\begin{align}
\gamma_{pp}\big(E_{n,s}\big)\geq K_{s+1}\-8(p-1)^2\sum_{b=1}^B (P_b+1)\cdot K_s-5B(p-1)^2-(p-1)\sum_{b=1}^B (P_b +1).\nonumber
\end{align}
\begin{flushright}
$\Box$
\end{flushright}
Based on this result,
\begin{lemma}\label{l:LowerBdd}
It holds that
$$\gamma_{pp}\big(\cap_{n\geq1}E_{n,s}\big)\geq K_{s+1}-8(p-1)^2\sum_{b=1}^B (P_b+1)\cdot K_s-5B(p-1)^2-(p-1)\sum_{b=1}^B (P_b +1).$$
Here, the set $E_{n,s}$ is defined as in (\ref{SetE}).
\end{lemma}
\noindent{\bf Proof.}
It is evident that
\begin{itemize}
  \item $E_{n+1,s}\subset E_{n,s}$;
  \item $E_{n,s}$'s are compact.
\end{itemize}
Then, $\cap_{n\geq1}E_{n,s}$ is non-empty and compact. Denote
$$g:=K_{s+1}-8(p-1)^2\sum_{b=1}^B (P_b+1)\cdot K_s-5B(p-1)^2-(p-1)\sum_{b=1}^B (P_b +1).$$
Let us argue by contradiction. Assume that $\gamma_{pp}(\cap_{n\geq1}E_{n,s})\leq g-1$.
Then, there is a $\varepsilon>0$ such that $\gamma_{pp}((\cap_{n\geq1}E_{n,s})_\varepsilon)\leq g-1$.

Now we claim that there is a $n_0>0$ such that $E_{n_0,s}\subset(\cap_{n\geq1}E_{n,s})_\varepsilon$. Otherwise, if for any $n\geq0$, $E_{n,s}\backslash (\cap_{n\geq1}E_{n,s})_\varepsilon\neq\emptyset$ and compact, then $\cap_{n\geq0}E_{n,s}\backslash (\cap_{n\geq1}E_{n,s})_\varepsilon\neq\emptyset$. This is a contradiction.

Hence, $g\leq\gamma_{pp}(E_{n_0,s})\leq\gamma_{pp}( (\cap_{n\geq1}E_{n,s})_\varepsilon)\leq g-1$. This is a contradiction again.
\begin{flushright}
$\Box$
\end{flushright}

Before proving Theorem \ref{t:existence}, let us summarize the properties we have known.
\begin{proposition}\label{prop:list}
It holds that
\begin{itemize}
  \item [$(a).$] The set $W_s:=\cap_{n\geq1}E_{n,s}\subset\partial\mathcal{A}\cap G_{K_{s+1}}$;
  \item [$(b).$] For any $U\in W_s$, any $t\geq0$, $n(\eta^t(U)_{(b-1)p+j})\equiv P_b$ for $b=1,\cdots,B$, $j=1,\cdots,p$ and $n(\eta^t(U)_{Bp+r})\equiv Q_r$ for $r=1,\cdots,R$. Moreover, for any $t\geq0$, each bump of each component of $\eta^t(U)$ has $L^4$-norm larger than $\rho$. Here, $\rho$ is the constant in Proposition \ref{p:L4};
  \item [$(c).$] For any $U\in W_s$, any $t>0$ and any $b=1,\cdots,B$ and any $j,j'=1,\cdots,p$ with $j\neq j'$, $4(P_{b}+1)\cdot K_s+2\leq n(\eta^t(U)_{(b-1)p+j}-\eta^t(U)_{(b-1)p+j})\leq 4(P_b+1) \cdot K_{s+1}+1$. Here, the upper bound is due to (\ref{ineq:ComparisonUpper});
  \item [$(d).$] $\gamma_{pp}(W_s)\geq K_{s+1}-8(p-1)^2\sum_{b=1}^B (P_b+1)\cdot K_s-5B(p-1)^2-(p-1)\sum_{b=1}^B (P_b +1)$.
\end{itemize}
\end{proposition}
\noindent{\bf Proof.}
Assertion $(1)$ holds because of the definition (\ref{SetE}). Assertion $(2)$ is due to $W_s\subset D_{K_{s+1}}$, (\ref{definitionDk}) and Proposition \ref{p:L4}. (\ref{SetE}) and Assertion $(a)$ imply Assertion $(3)$. Lemma \ref{l:LowerBdd} ensures Assertion $(d)$.

\begin{flushright}
$\Box$
\end{flushright}

\subsection{Proof of Theorem \ref{t:existence}}\label{Subsection:proof}

To prove the first part of Theorem \ref{t:existence}, it is sufficient to notice that
\begin{align}
K_{s+1}-8\sum_{b=1}^B (P_b+1)\cdot K_s(p-1)^2-5B(p-1)^2-\sum_{b=1}^B (P_b +1)(p-1)=1\nonumber
\end{align}
with $K_1=8\sum_{b=1}^B (P_b +1)p(p+1)+5Bp^2$. This is due to (\ref{def:Ks}). Therefore, we get $W_s\neq\emptyset$ for any $s$.
In order to continue the discussion, let us introduce the following claim concerning the relation between the $\omega$-set and the equilibrium points of Problem (\ref{e:BBBB}).
\begin{claim}\label{c:equilibrium}
There exists a $U_\infty\in\omega(W_s)$ solving Problem (\ref{e:AAAA}). Here, the set $W_s$ is defined in Proposition \ref{prop:list}.
\end{claim}
\noindent{\bf Proof of Claim \ref{c:equilibrium}.}
We apply an idea in \cite{CMT2000,LiWang2021}. Since $\gamma_{pp}(W_s)\geq 1$ due to (\ref{ASSUMPTION2}) and Assertion $(d)$ of Proposition \ref{prop:list} in the last subsection, $W_s\neq\emptyset$. For any $U_0\in W_s$, it is evident that $U_0\in\partial\mathcal{A}\cap(H_{0,r}^2(\Omega))$. Therefore, $T(U_0)=\infty$ and $\inf_{t\geq 0}I(\eta^t(U_0))>0$. Proposition \ref{p:L2partial_t} implies that
\begin{align}
\sum_{j=1}^N\int_0^\infty|\partial_t \eta^t(U_0)_j|_2^2dt=I(U_0)-\lim_{t\to+\infty}I(\eta^t(U_0))<+\infty. \nonumber
\end{align}
Here, $\eta^t(U_0)$ is the $j$-th component of the vector-valued function $\eta^t(U)$. Then, there exists a sequence $\{t_n\}_n\subset\mathbb{R}$ with $|\partial_t\eta^{t_n}(U_0)|_2\to0$ as $n\to\infty$. Therefore, for the same sequence $\{t_n\}$, it holds that
\begin{align}
\nabla I(\eta^{t_n}(U_0))\to0\mbox{ in }(H_r^{-1}(\Omega))^N.\nonumber
\end{align}
It follows that $\{\eta^{t_n}(U_0)\}_n$ is a $(PS)$ sequence. Since the functional $I$ satisfies $(PS)$ condition obviously, there exists a function $U_\infty\in (H_{0,r}^1(\Omega))^N$ such that $\eta^{t_n}(U_0)\to U$ in $(H_{0,r}^1(\Omega))^N$ and $U_\infty$ solves Problem (\ref{e:AAAA}). Using the definition of omega set, $U_\infty\in \omega(W_s)$.

\begin{flushright}
$\Box$
\end{flushright}

\begin{claim}\label{c:nodes}
For the solution $U_\infty$ to Problem (\ref{e:AAAA}) in Claim \ref{c:equilibrium}, it holds that
\begin{itemize}
  \item [$(1).$] for any $b=1,\cdots,B$ and any $i=1,\cdots,p$, $n((U_\infty)_{(b-1)p+i})=P_b$;
  \item [$(2).$] for any $r=1,\cdots,R$, $n((U_\infty)_{Bp+r})=Q_r$.
\end{itemize}
\end{claim}
\noindent{\bf Proof of Claim \ref{c:nodes}.}
We follow the notation of Claim \ref{c:equilibrium}.
Using Proposition \ref{prop:list}, it is known that for any $n\geq0$, $n(\eta^{t_n}(U)_{(b-1)p+j})\equiv P_b$ for $b=1,\cdots,B$, $j=1,\cdots,p$ and $n(\eta^{t_n}(U)_{Bp+r})\equiv Q_r$ for $r=1,\cdots,R$. Moreover, each bump of each component of $\eta^{t_n}(U)$ has $L^4$-norm smaller than $\rho$. Here, $\rho$ is the constant in Proposition \ref{p:L4}.

Recall that $\eta^{t_n}(U_0)\to U_\infty$ in $(H_{0,r}^1(\Omega))^N$.
If there exists a $b_0=1,\cdots,B$ and $j_0=1,\cdots,p$ such that $n((U_\infty)_{(b_0 -1)p+j_0})<n(\eta^{t_n}(U)_{(b-1)p+j})=P_b$, there exists a $r_0>0$ such that $\partial B_{r_0}(0)\subset\overline{\Omega}$ and $(U_\infty)_{(b_0 -1)p+j_0} |_{|x|=r_0}=|\nabla (U_\infty)_{(b_0 -1)p+j_0}|_{|x|=r_0}=0$. Using the unique solvability of ODE, $(U_\infty)_{(b_0 -1)p+j_0}\equiv 0$ in $\Omega$. This contradicts with the construction of $W_s$. A similar argument can be proceed for $n((U_\infty)_{Bp+r})=Q_r$ with $r=1,\cdots,R$. Moreover, due to the $H^1$-convergence of $\eta^{t_n}(U_0)$, it is evident that each bump of each component of $U_\infty$ has $L^4$-norm greater than $\frac{\rho}{2}$

\begin{flushright}
$\Box$
\end{flushright}

\begin{remark}
Claim \ref{c:nodes} proves Assertions $(1)$ and $(2)$ of Theorem \ref{t:existence}.
\end{remark}

In the next step, we estimate a part of the components of $U_\infty$. To be precise,
\begin{claim}\label{c:compare0}
For the solution $U_\infty$ to Problem (\ref{e:AAAA}) in Claim \ref{c:equilibrium}, for any $b=1,\cdots,B$ and any $j,j'=1,\cdots,p$ with $j\neq j'$, $4(P_{b}+1)\cdot K_s+2\leq n((U_\infty)_{(b-1)p+j}-(U_\infty)_{(b-1)p+j'})\leq 4(P_b+1) \cdot K_{s+1}+1$.
\end{claim}

This can be proved via a similar argument as in Claim \ref{c:nodes}.
Moreover, we have
\begin{claim}\label{c:compare}
It holds that
\begin{itemize}
  \item [$(1).$] for any $b_1,b_2=1,\cdots,B$ and any $i_1,i_2=1\cdots,p$ with $b_1\neq b_2$, $n((U_\infty)_{(b_1 -1) p+i_1}-(U_\infty)_{(b_2 -1) p+i_2})\leq P_{b_1}+P_{b_2}+1$;
  \item [$(2).$] for any $b=1,\cdots,B$, $i=1,\cdots,p$ and $r=1,\cdots,R$, $n((U_\infty)_{(b-1)p+i}-(U_\infty)_{Bp+r})\leq P_b +Q_r +1$;
  \item [$(3).$] for any $r_1 ,r_2=1,\cdots,R$ with $r_1\neq r_2$, $n((U_\infty)_{Bp+r_1} - (U_\infty)_{Bp+r_2})\leq Q_{r_1}+Q_{r_2}+1$.
\end{itemize}
\end{claim}
\noindent{\bf Proof of Claim \ref{c:compare}.}
We follow the notation of Claim \ref{c:equilibrium}.
Due to Corollary \ref{c:nodalnumber}, for any $b_1,b_2=1,\cdots,B$ and any $i_1,i_2=1,\cdots,p$
$$n((U_s)_{pb_1 +i_1-p}-(U_s)_{(b_2 -1)p +i_2})\leq n((U_\infty)_{(b_1 -1)p +i_1}-(U_\infty)_{(b_2 -1)p +i_2})\leq P_{b_1} +P_{b_2} +1.$$
For any $r,r'=1,\cdots,R$ with $r\neq r'$,
$$n((U_s)_{r}-(U_s)_{r'})\leq n((U_\infty)_{r}-(U_\infty)_{r'})\leq Q_r +Q_{r'} +1.$$
For any $b=1,\cdots,B$, $i=1,\cdots,p$ and $r=1,\cdots,R$,
$$n((U_s)_{(b-1)p+i}-(U_s)_{Bp+r})\leq n((U_\infty)_{(b-1)p+i}-(U_\infty)_{Bp+r})\leq P_b +Q_r +1.$$
\begin{flushright}
$\Box$
\end{flushright}

\noindent{\bf Proof of Theorem \ref{t:existence}.}
Theorem \ref{t:existence} follows from Claims \ref{c:equilibrium}, \ref{c:nodes}, \ref{c:compare0} and \ref{c:compare}.
\begin{flushright}
$\Box$
\end{flushright}

\section{Proof of Theorem \ref{t:nonexistence}}

In this section, we prove Theorem \ref{t:nonexistence}.

\noindent{\bf Proof of Theorem \ref{t:nonexistence}.}
We first verify Assertion $(1)$ of Theorem \ref{t:nonexistence}.
Consider
\begin{equation}\label{e:CCC}
    \left\{
   \begin{array}{lr}
     -{\Delta}u_{j}+ u_{j}= u^{3}_{j}+\sum_{i=1, i\neq j}^N\beta u_{j}u_{i}^{2} \,\,\,\,\,\,\,  \mbox{in}\ \Omega ,\\
     0<u_{j}\in H_{0}^{1}(\Omega), \,\,\,\,\,\,\,\,j=1,\dots,N.
   \end{array}
   \right.
\end{equation}
Here, the domain $\Omega\subset\mathbb{R}^n$ for $n=2,3$ is any domain with smooth boundary. The constant satisfies $\beta\leq-1$.
\begin{claim}\label{c:nonexistence}
For a solution $(u_1,\cdots,u_N)$ to Problem (\ref{e:CCC}), $w_{ij}=u_i - u_j$ must change its sign for $i,j=1,\cdots,N$ and $i\neq j$.
\end{claim}
Otherwise, without loss of generality, let $u_1\geq u_2>0$ in $\Omega$. Then,
\begin{align}
\int|\nabla u_2|^2+\int u^2_2 & =\int u_2^4+\beta\int u_1^2 u_2^2+\beta\sum_{i\geq3}\int u_i^2 u_1^2\nonumber\\
&\leq(1+\beta)\int u_2^4\leq0.\nonumber
\end{align}
This is a contradiction.

Now we check Assertion $(2)$ of Theorem \ref{t:nonexistence}.
Without loss of generality, let us consider $u_1$ and $u_2$. Suppose that $u_2(0)>u_1(0)>0$ and $z_1<z_2<\cdots<z_{n(u_1)}$ are the zeroes of $u_1$.

Firstly, notice that on the interval $(0,z_1)$, the graphs of $u_1$ and $u_2$ must intersect. Otherwise, on $(0,z_1)$
\begin{itemize}
  \item $0<u_1<u_2$;
  \item $-\Delta u_1 + u_1= u_1^3+\beta\sum_{i\geq2}u^2_i u_1$.
\end{itemize}
Then, $\int_{B(0,z_1)}|\nabla u_1|^2+ u_1^2\leq(1 +\beta)\int_{B(0,z_1)}u_1^2 u_2^2\leq0$. For $k=1,\cdots,\big[\frac{n(u_1)-1}{2}\big]$, on the interval $[z_{2k-1},z_{2k+1})$, we can proceed a similar argument and obtain an intersection of $u_1$ and $u_2$. This proves the theorem.

\begin{flushright}
$\Box$
\end{flushright}

\appendix
\section{Proof of (\ref{ineq:genus001}), (\ref{ineq:genus11}) and (\ref{ineq:genus22})}
In this appendix, we prove (\ref{ineq:genus001}), (\ref{ineq:genus11}) and (\ref{ineq:genus22}).

\vs

\noindent{\bf Proof of (\ref{ineq:genus001}).}
In this part, we prove that
\begin{align}
\gamma_{\sigma_p}(\partial\mathcal{O})\leq K\Big(\sum_{b=1}^{B}P_b +\sum_{r=1}^R Q_r +B+R\Big).\nonumber
\end{align}
Here,
$\mathcal{O}=\{z\in X_1| |h(z)_{j,q}|_4\leq{2\varepsilon}\mbox{ for any admissible }(j,q)\}$ and
\begin{align}
X_1&=\Big\{z=(z^{(1)},z^{(2)})\in\mathbb{C}^{K(p\sum_{b=1}^B P_b+\sum_{r=1}^R Q_r+N)}\Big|\mbox{for any }b=1,\cdots,B;\,q=1,\cdots,P_b+1;\nonumber\\
&\quad k=1,\cdots,K,\mbox{ we have }\mbox{arc}(z_{(b-1)p+1,q,k}^{(1)})=\mbox{arc}(z_{(b-1)p+2,q,k}^{(1)}) +\frac{2\pi}{p}=\cdots\nonumber\\
&\quad=\mbox{arc}(z^{(1)}_{bp,q,k})+\frac{2\pi(p-1)}{p}\Big\}.\nonumber
\end{align}
The numbers $B$, $P_b$, $R$ and $Q_r$ are prescribed in Theorem \ref{t:existence}.
It is sufficient to construct a mapping
\begin{align}
\Phi:X_1\to\mathbb{C}^{K(\sum_{b=1}^{B}P_b +\sum_{r=1}^R Q_r +B+R)}\nonumber
\end{align}
with
\begin{itemize}
  \item [$(a).$] $\Phi(e^{\frac{2\pi i}{p}}z)=e^{\frac{2\pi i}{p}}\Phi(z)$;
  \item [$(b).$] $\Phi^{-1}(0)=0$.
\end{itemize}

To do this, we only need to define
\begin{align}
\Phi(z_{(b-1)p+j,q,k}^{(1)},z_{Bp+r,q,k}^{(2)})=(\sum_{j=0}^{p-1}e^{\frac{2\pi i}{p}\cdot j}z_{(b-1)p+j,q,k}^{(1)},z_{Bp+r,q,k}^{(2)}).\nonumber
\end{align}
This completes the proof of (\ref{ineq:genus001}).

\begin{flushright}
$\Box$
\end{flushright}

Before we prove (\ref{ineq:genus11}) and (\ref{ineq:genus22}), let us briefly recall the notations.

\begin{itemize}
  \item $R_{(b-1)p+j,q}^{2\varepsilon}=\{z\in X_1||h(z)_{(b-1)p+j,q}|_4\leq{2\varepsilon}\}$ for $q=1,\dots,P_{b}+1$ and $j=1,\dots,N$;
  \item $R_{Bp+r,q}^{2\varepsilon}=\{z\in X_1||h(z)_{Bp+r,q}|_4\leq{2\varepsilon}\}$ for $r=1,\cdots,R$ and $q=1,\cdots,Q_r$. The mapping $h$ is defined as in (\ref{def:h});
  \item $A_1=\{z\in X_1|\exists b\in\{1,\cdots,B\},\exists q\in\{1,\cdots, P_{b}+1\},\exists i,i'\in\{1,\cdots,p\}\mbox{ with }i\neq i'\mbox{ such that }z\in R_{(b-1)p + i, q}^{2\varepsilon}\cap S_{(b-1)p +i', q}^{2\varepsilon}\}$;
  \item $A_2=\cup_{b=1}^B\cup_{q=1}^{P_b +1}\cap_{j=1}^p R_{(b-1)p + j,q}^{2\varepsilon}$;
  \item $A_3=\cup_{r=1}^R\cup_{q=1}^{Q_r + 1} R_{Bp+r,q}^{2\varepsilon}$.
  \item $A_2=\cup_{b=1}^B\cup_{q=1}^{P_b +1}\cap_{j=1}^p R_{(b-1)p + j,q}^{2\varepsilon}$;
  \item $A_3=\cup_{r=1}^R\cup_{q=1}^{Q_r + 1} R_{Bp+r,q}^{2\varepsilon}$.
\end{itemize}

\vs

\noindent{\bf Proof of (\ref{ineq:genus11}).}
Recall that for any $z\in A_1$, there are $b_U=1,\cdots,B$, $q_U=1,\cdots, P_{b_U}+1$ and $i_U,i'_{U}=1,\cdots,p$ with $i_U\neq i'_U$ such that $z\in R_{(b_U -1)p + i_U , q_U}^{2\varepsilon}\cap S_{(b_U -1)p +i'_U , q_U}^{2\varepsilon}$. For any $b=1,\cdots,B$, $q=1,\cdots,P_b+1$ and $j=1,\cdots,p$, let us define
\begin{align}
A_{1,(b-1)p+j,q}=A_1\cap R_{(b-1)p+j,q}^{2\varepsilon-\delta}.\nonumber
\end{align}
It is evident that
\begin{align}
A_1=\cup_{b=1}^B\cup_{q=1}^{P_b+1}\cup_{j=0}^{p-1}\sigma_p^j(A_{1,(b-1)p+1,q})\nonumber
\end{align}
and
\begin{align}
\cap_{j=0}^{p-1}\sigma_p^j(A_{1,(b-1)p+1,q})=\emptyset.\nonumber
\end{align}
Using Assertion $(7)$ of Proposition \ref{prop:genus}, we get
\begin{align}
\gamma_{\sigma_p}\big(\cup_{j=0}^{p-1}\sigma_p^j(A_{1,(b-1)p+1,q})\big)\leq p-1.\nonumber
\end{align}
Therefore, by Assertion $(2)$ of Proposition \ref{prop:genus},
\begin{align}
\gamma_{\sigma_p}(A_1)\leq (p-1)\Big(\sum_{b=1}^B P_b+B\Big).\nonumber
\end{align}

\begin{flushright}
$\Box$
\end{flushright}

\vs

\noindent{\bf Proof of (\ref{ineq:genus22}).}
In this part, we make use of an idea in \cite{LiuLiuWang2015}. Recall that
\begin{itemize}
  \item $A_2=\cup_{b=1}^B\cup_{q=1}^{P_b +1}\cap_{j=1}^p R_{(b-1)p + j,q}^{2\varepsilon}$;
  \item $A_3=\cup_{r=1}^R\cup_{q=1}^{Q_r + 1} R_{Bp+r,q}^{2\varepsilon}$.
\end{itemize}
We want to check that
\begin{align}
\gamma_{\sigma_p}\big((A_2 \cup A_3)\cap\partial\mathcal{O}\backslash A_1\big)\leq K\Big(\sum_{b=1}^{B}P_b +\sum_{r=1}^R Q_r +B+R-1\Big).\nonumber
\end{align}

To begin with, we introduce a new index set.
\begin{align}
\mathcal{I}=\{&(1,1),\cdots,(1,P_1+1);\cdots;(B,1),\cdots,(B,P_B+1);\nonumber\\
&(B+1,1),\cdots,(B+1,Q_1+1);\cdots;(B+R,1),\cdots,(B+R,Q_R+1)\}.\nonumber
\end{align}
For any $(j,p)\in\mathcal{I}$, define
\begin{align}
T_{j,q}=\cap_{j=1}^p R_{(b-1)p+j,q}^{2\varepsilon}\nonumber
\end{align}
if $j=1,\cdots,B$ and
\begin{align}
T_{j,q}=R_{B(p-1)+j,q}^{2\varepsilon}\nonumber
\end{align}
of $j=B+1,\cdots,B+R$. Then,
\begin{align}
(A_2\cup A_3)\cap\partial\mathcal{O}\backslash A_1\subset\cup_{l=1}^{\sum_{b=1}^N P_b+\sum_{r=1}^R Q_r+B+R-1}B_l\nonumber
\end{align}
with
\begin{itemize}
  \item $B_l=\cup_{s\in S_l}\Big[(\cap_{(j,q)\in s}\partial T_{j,q})\cap(\cap_{(j,q)\in s^c}T_{j,q})\Big]$;
  \item $S_l=\{s\subset\mathcal{I}|\# s=l\}$;
  \item for any $s\in S_l$, $s^c:=\mathcal{I}\backslash s$.
\end{itemize}
Notice that
\begin{align}
B_{\sum_{b=1}^N P_b+\sum_{r=1}^R Q_r+B+R}\subset A_1\cup L.\nonumber
\end{align}
For any $l=1,\cdots,\sum_{b=1}^N P_b+\sum_{r=1}^R Q_r+B+R-1$, define the mapping
\begin{align}
f_l:B_l\to\mathbb{C}^K\nonumber
\end{align}
by
\begin{align}
f_l(z)=\sum_{s\in S_l}(\Phi(z))_{i(s)}\cdot d\Big(z,\partial(\cap_{(j,q)\in s^c}T_{j,q})\Big).\nonumber
\end{align}
Here, $i(s)$ is the first couple in $s$ by the dictionary order and
\begin{align}
(\Phi(z))_{(j,q)}=\big((\Phi(z))_{q,1},\cdots,(\Phi(z))_{q,K}^j\big).\nonumber
\end{align}
It is evident that
\begin{claim}
\begin{itemize}
  \item [$(1)$] $f_l$ is $\sigma_p$-equiv-variant;
  \item [$(2)$] $f_l\neq0$ on $B_l$.
\end{itemize}
\end{claim}
This implies that
\begin{corollary}
$\gamma_{\sigma_p}((A_2\cup A_3)\cap\backslash A_1)\leq K\Big(\sum_{b=1}^{B}P_b +\sum_{r=1}^R Q_r +B+R-1\Big)$.
\end{corollary}
This completes the proof.

\begin{flushright}
$\Box$
\end{flushright}

\noindent{\bf Acknowledgement.}
Li is supported by FAPESP Proc 2022/15812-0.

\end{sloppypar}

\end{document}